\documentclass[12pt]{amsart}
\usepackage{amsmath,amsthm,amsfonts,amssymb,eucal}

\renewcommand {\a}{ \alpha }
\renewcommand{\b}{\beta}
\newcommand{\e}{\epsilon}
\newcommand{\vare}{\varepsilon}
\newcommand{\g}{\gamma}
\newcommand{\G}{\Gamma}

\newcommand{\vark}{\varkappa}
\renewcommand{\d}{\delta}
\newcommand{\s}{\sigma}
\renewcommand{\l}{\lambda}

\newcommand{\z}{\zeta}
\renewcommand{\t}{\theta}

\newcommand{\T}{\Theta}

\newcommand{\om}{\omega}

\newcommand{\R}{ \mathbb R}

\newcommand {\gb}{\mathfrak b}
\newcommand{\gl}{\mathfrak l}

\newcommand {\ba}{\mathbf a}

\newcommand {\BA}{\mathbf A}

\newcommand{\CV}{\mathcal V}

\newcommand{\CB}{\mathcal B}

\newcommand{\CP}{\mathcal P}

\newcommand{\CA}{\mathcal A}

\newcommand{\CC}{\mathcal C}

\newcommand{\CE}{\mathcal E}

\newcommand{\plainC}[1]{\textup{{\textsf{C}}}^{#1}}

\newcommand{\plainH}[1]{\textup{{\textsf{H}}}^{#1}}

\newcommand{\plainL}[1]{\textup{{\textsf{L}}}^{#1}}

\DeclareMathOperator{\modulo}{{mod}}

 \DeclareMathOperator{\dc}{d}
 
\DeclareMathOperator{\sign}{\boldsymbol{\mathfrak{s}}}

\newtheorem{thm}{Theorem}[section]
\newtheorem{cor}[thm]{Corollary}
\newtheorem{lem}[thm]{Lemma}
\newtheorem{prop}[thm]{Proposition}
\newtheorem{cond}[thm]{Condition}

\theoremstyle{definition}
\newtheorem{defn}[thm]{Definition}

\theoremstyle{remark}
\newtheorem{rem}[thm]{Remark}

\numberwithin{equation}{section}

\newcommand{\thmref}[1]{Theorem~\ref{#1}}

\newcommand{\lemref}[1]{Lemma~\ref{#1}}

\newcommand{\bee}{\begin{equation}}
\newcommand{\ene}{\end{equation}}
\newcommand{\bes}{\begin{split}}
\newcommand{\ens}{\end{split}}

\newcommand{\bet}{\begin{tm}}
\newcommand{\ent}{\end{tm}}
\newcommand{\bel}{\begin{lm}}
\newcommand{\enl}{\end{lm}}
\newcommand{\bec}{\begin{cor}}
\newcommand{\enc}{\end{cor}}
\newcommand{\bep}{\begin{pr}}
\newcommand{\enp}{\end{pr}}
\newcommand{\ber}{\begin{rem}}
\newcommand{\enr}{\end{rem}}

\newcommand{\Z}{\mathbb Z}

\newcommand{\res}{\!\restriction\!}
\newcommand{\D}{\Delta}
\newcommand{\CF}{\mathcal F}

\DeclareMathOperator{\dom}{{Dom}}

\DeclareMathOperator{\gen}{{Gen}}

\makeatletter
\def\square{\RIfM@\bgroup\else$\bgroup\aftergroup$\fi
  \vcenter{\hrule\hbox{\vrule\@height.6em\kern.6em\vrule}\hrule}\egroup}
\makeatother

\begin{document}
\hoffset -4pc
\title[Schr\"odinger operator on trees]
{Schr\"odinger operator on homogeneous metric trees: spectrum in
gaps}
\author[A. V. Sobolev]
{Alexander V. Sobolev}
\address{Centre for Mathematical
Analysis and Its Applications\\ University of Sussex\\ Falmer,
Brighton\\ BN1 9QH, UK} \email{A.V.Sobolev@sussex.ac.uk}
\author[M. Solomyak]
{Michael Solomyak}
\address
{Department of Mathematics
\\ Weizmann Institute\\ Rehovot\\ Israel}
\email{solom@wisdom.weizmann.ac.il}
\thanks{}
\date{31 July 2001}
\begin{abstract}
The paper studies the spectral properties of the Schr\"odinger
operator $\BA_{gV} = \BA_0 + gV$ on a homogeneous rooted metric
tree, with a decaying real-valued potential $V$ and a coupling
constant $g\ge 0$. The spectrum of the free Laplacian $\BA_0 =
-\boldsymbol{\D}$ has a band-gap structure with a single
eigenvalue of infinite multiplicity in the middle of each finite
gap. The perturbation $gV$ gives rise to extra eigenvalues in the
gaps. These eigenvalues are monotone functions of $g$ if the
potential $V$ has a fixed sign. Assuming that the latter condition
is satisfied and that $V$ is symmetric, i.e. depends on the
distance to the root of the tree, we carry out a detailed
asymptotic analysis of the counting function of the discrete
eigenvalues in the limit $g\to\infty$. Depending on the sign and
decay of $V$, this asymptotics is either of the Weyl type or is
completely determined by the behaviour of $V$ at infinity.
\end{abstract}

\maketitle
\section{Introduction}

Counting the number of eigenvalues of a perturbed operator,
appearing in the spectral gaps of the unperturbed one, is a
classical problem. It was extensively investigated both in the
general operator-theoretic setting \cite{BirAdv} and in
applications to various specific problems of Mathematical Physics
(the Hill operator, \cite{DH}, \cite{Sob}; the Dirac operator,
\cite{Klaus}, \cite{BirLap}; the  periodic Schr\"odinger and
magnetic Schr\"odinger operators, \cite{ADH}, \cite{BirRai};
waveguide-type operators, \cite{BirSol}, etc.) In this paper we
study a new problem of this type, which only recently attracted
the attention of specialists: our unperturbed operator is the
Laplacian on a homogeneous rooted metric tree $\G$. In general, a
metric tree is a tree whose edges are viewed as non-degenerate
line segments, rather than pairs of vertices, as in the case of
the standard ( combinatorial ) trees. This difference is reflected
in the nature of the corresponding Laplacian. For a combinatorial
tree this is the discrete Laplacian, whereas the Laplacian $\BA_0
= -\boldsymbol\D$ on a metric tree is represented by a family of
the operators $-d^2/dx^2$ on its edges, complemented by the
Kirchhoff matching conditions at the vertices. The Laplacian on
the homogeneous metric tree has very specific spectral properties
which we describe later on in details. In particular, the spectrum
has the band-gap structure, with a single eigenvalue of infinite
multiplicity in each finite gap. For some other operators on a
homogeneous tree, having similar nature, the band-gap structure of
the spectrum was established earlier by R.Carlson \cite{carlson1}.
In the present paper we study the properties of the perturbed
operator $\BA_{gV} = \BA_0 + gV$ where $V$ is a decaying
real-valued potential, and $g\ge 0$ is a coupling constant. The
potential $V$ is assumed to be symmetric, i.e. dependent only on
the distance $|x|$ between $x\in\G$ and the root of $\G$. This
perturbation may produce extra eigenvalues in the gaps of $\BA_0$.

For an ``observation point''
$\l$ inside a gap we denote by $M(\l; \BA_{gV})$
the number of the eigenvalues
of the operator $\BA_{\a V}$,
crossing $\l$ as $\a$ varies from $0$
to $g$. For any two points $\l_1,\l_2,\ \l_1<\l_2$
lying in the same gap,
we denote by $N(\l_1,\l_2; \BA_{gV})$
he number of eigenvalues of this
operator on the interval $(\l_1,\l_2)$;
see Sect. \ref{tree4:sect}
for more precise definitions. We are
interested in the limiting behaviour
of these quantities as the coupling constant $g$
tends to infinity. Compared
with other problems of this type,
mentioned in the beginning of the
Introduction, this problem has
 many new important features.

\vskip0.3cm

The starting point of our investigation is a direct
decomposition of the Sobolev space $\plainH{1,0}(\G)$
on the homogeneous tree $\G$. This decomposition is
orthogonal with respect
to the inner products $\int_\G u'\overline{v'} dx$ and
$\int_\G V(x)u\overline v dx$ for all
symmetric weight functions $V$
simultaneously.
Therefore, it reduces the Laplacian and any
Schr\"odinger operator $\BA_V$ with a
symmetric potential $V$. This
decomposition was constructed in the paper
\cite{NS} and has proved very
useful for spectral theory of this
class of operators. Later it
was re-discovered by
R.Carlson \cite{carlson2} in a somewhat different setting.

The parts of the operator $\BA_{V}$ in each component of
the said orthogonal decomposition
turn out to be unitary equivalent to
second order differential
operators $A_{V_k}$, $k=0,1,\ldots$
of the Sturm-Liouville type in the space
$\plainL2(\R_+)$. The potentials $V_k$ are obtained
from the original potential $V$
by ``shifting'' the variable:
$V_k(t)=V(t+k)$, $k=0,1,\ldots$.
The operators $A_{V_k}$ act as
$ -d^2/dx^2 + V_k$, but in contrast to the standard
Sturm-Liouville problem, the
description of the operator
domain of $A_{V_k}$ involves specific
matching conditions at the points
$t_n=n$, $n=1,2\ldots$ (see Section \ref{tree2:sect}).
Each component  $A_{V_k}$, $k\ge 1$,
enters  $\BA_{V}$
with the multiplicity
\begin{equation*}
n_0 = 1,\ \ n_k = b^k-b^{k-1},\ k\ge 1.
\end{equation*}
Here $b\ge2$ is the integer-valued parameter
(the branching number)
which characterizes the homogeneous
tree completely, see the definition in
Subsect. \ref{tree2.1:subsect}.

The above orthogonal
decomposition plays a central
role in our approach. First of all, it
allows us to calculate  the spectrum of $\BA_0$ explicitly
(see Theorem \ref{sp:thmA1}):
it consists of the bands
$\bigl[(\pi(l-1)+\t)^2,(\pi l-\t)^2\bigr]$,\
$\t=\arccos\bigl(2(b^{1/2}+b^{-1/2})^{-1}\bigr)$,  and
the eigenvalues $\l_l=(\pi l)^2$, $l\in \mathbb N$.
Besides, for $V=0$ all the components $A_{V_k}=A_0$
are identical, so that the spectrum
is of infinite multiplicity.
For the perturbed operator this
decomposition leads to the representation
\begin{equation}\label{yes:eq}
M(\l; \BA_{gV}) = \sum_{k \ge 0} n_k M(\l; A_{gV_k}).
\end{equation}
A similar formula also holds for
$N(\l_1,\l_2;\BA_{gV})$.
The presence of the exponentially growing
factors $n_k$ hampers the study of
these sums.  Remembering that the numbers $n_k$
reflect the geometry of the
tree, rather than the properties
of the potential $V$, we also
study the counting functions $\widetilde M$ and $\widetilde N$
ignoring the exponential
multiplicities $n_k$. Precisely, we introduce
\begin{equation}\label{no:eq}
\widetilde M(\l; \BA_{gV}) = \sum_{k\ge 0} M(\l; A_{gV_k}),
\end{equation}
and the quantity $\widetilde N$ defined in a similar way.

Clearly, the study of the four
functions $M,N,\widetilde M,\widetilde N$
reduces to that of the individual counting functions
$M(\l; A_{gV_k}), N(\l_1, \l_2; A_{gV_k})$
for the operators $A_{V_k}$.
A similar problem for the
classical Hill operator was investigated in
\cite{DH}
and, in a more detailed way, in \cite{Sob}.
The general strategy adopted in
\cite{Sob} applies to the operators
$A_{gV_k}$ with only minor changes.
However, here a new problem emerges:
in order  to obtain the asymptotic formulas for the sums
\eqref{yes:eq}, \eqref{no:eq} one needs
an asymptotics of $M(\l; A_{gV_k})$ jointly
in two parameters: $g$ and $k$, with a good control of the
remainder estimate.
In solving this new
problem we see the main technical novelty of the paper.

In the paper we obtain several
results of rather different type.
In this introduction we do not describe them
in detail, but concentrate on their principal features.
More extended comments are given in the main text.
Also, for the sake of discussion we restrict
ourselves to the functions $M(\l; \BA_{V})$
and $\widetilde M(\l; \BA_{gV})$ only.

First of all, as in the case of the ``classical''
Hill operator problem (see \cite{Sob}),
we observe that the behaviour of $M, \widetilde M$
is in general radically different for non-positive and
non-negative potentials. More precisely, if $V\le 0$ decays
sufficiently quickly at infinity, then the
asymptotics is governed by an appropriate Weyl-type
formula, and thus it depends on the values of $V(t)$ at all
points $t\in \R$, and contains no information
on the spectrum of the unperturbed operator.
On the contrary, for a non-negative $V$
the asymptotics is determined by
the fall-off of $V$ at infinity,
and as a rule, depends heavily on some spectral characteristic
of the operator $\BA_0$. For instance,
the behaviour of $\widetilde M(\l)$ for $V\ge 0$ is
described by an integral of
the density of states for $A_0$. Similar type of asymptotics
is also observed for the potentials $V\le 0$ whose decay at infinity is
slow in some specified sense.

In accordance with this general observation our
study of the asymptotics is
divided in several parts. We begin in Sect. \ref{tree4:sect}
by specifying the conditions on a non-positive
 symmetric potential $V$
 that guarantee the validity of the Weyl-type
asymptotics. Further on,
we proceed to the cases when the Weyl formula fails and
the asymptotics is determined by the
behaviour of $V$ at infinity. Here we investigate two
types of potentials: power-like and exponentially decaying.
In Sect. \ref{tree5:sect} we state the results for the functions
$\widetilde M, \widetilde N$.
A common feature of the asymptotic formulae
in Sect. \ref{tree5:sect}
is that virtually all of them contain the
density of states for the operator $A_0$.
It is also worth pointing out that
the power-like potentials induce
a power-like growth of $\widetilde M$ as $g\to\infty$,
whereas the exponential potentials give rise to a logarithmic growth.

The study of the sum \eqref{yes:eq}
is postponed until Sect. \ref{tree9:sect} as it
calls for different techniques
and is less complete.
For the power-like potentials
we are able to establish the asymptotics only for
the quantity
$\ln M(\l)$. For the exponential potentials
we provide more detailed asymptotic information. This is
possible due to the ``self-similarity''
of the exponential function. This property allows us to
rewrite the formula \eqref{yes:eq}
for the function $M(\l, \BA_{gV})$ in a form which
can be interpreted as a Renewal Equation (see \cite{F}, \cite{LV}).
Then the Renewal Theorem ensures a specific asymptotic behaviour of
$M(\l)$.

Let us briefly outline the contents of the remaining sections.
In Sect. \ref{tree2:sect} we describe
the basic orthogonal
decomposition of the space $\plainH{1,0}(\G)$ and
also the parts of the operator
$\BA_V$  in its components.
In Sect. \ref{tree3:sect} we
calculate the spectrum of the
Laplacian on $\G$. Here we also carry out a
detailed analysis of the density of states for the operator
$A_0$. This function is involved in the asymptotic
formulae for $\widetilde M$
in the non-Weyl situation.
As was mentioned earlier, the study of the perturbed operator
$\BA_V$ starts in Sect. \ref{tree4:sect}
where the Weyl's asymptotics is established.
The main results on the
non-Weyl asymptotics for the functions
$\widetilde M$ and $\widetilde N$ are collected in
Sect. \ref{tree5:sect}.
Their proofs are given in Sect. \ref{tree8:sect},
preceded by necessary technical
preliminaries in Sect. \ref{tree6:sect},
\ref{tree7:sect}.
The last Sect. \ref{tree9:sect} is devoted to the analysis
of the functions $M$, $N$.

\section{Laplacian on a homogeneous tree and its decomposition}
\label{tree2:sect}

\subsection{ Homogeneous trees and Laplacians on them}
\label{tree2.1:subsect}
Let $\G$ be a rooted tree with the root $o$,
the set of vertices $\CV(\G)$
and the set of edges $\CE(\G)$.
We suppose that the length of each edge
$e$ is equal to $1$. Given two points
$y,z\in\G$, we
write $y\preceq z$ if $y$ lies on
the unique simple path connecting $o$ with $z$;
let $|z|$ stand for the length
of this path. We write $y\prec z$ if
$y\preceq z$ and
$y\neq z$. The
relation $\prec$ defines on $\G$ a partial ordering.
If $y\prec z$,
we denote
\begin{equation*}
\langle y,z\rangle := \{ x\in\G:\;y\preceq x\preceq z \}.
\end{equation*}
In particular, if $e=\langle v,w\rangle$ is an edge, we call
$v$ its initial point
and say that $e$ emanates from $v$ and terminates at $w$.

For any $v\in\CV(\G)$
the number $|v|$ is a non-negative integer;
we call it \textsl{generation}
of $v$ and denote $\gen(v)$.
For an edge $e\in\CE(\G)$ $\gen(e)$
is defined as the generation of its initial point.

Let an integer $b>1$ be given.
We suppose that for each vertex $v\neq o$
there are exactly
$b$ edges emanating from $v$.
We denote them $e_v^1,\dots,e_v^b$ and write
$e_v^-$ for the edge terminating at $v$. We call $b$
\textsl{the branching number} of $\G$.
We always suppose that only one
edge emanates from
the root $o$.
Thus, the tree $\G$ is fully
determined
by the parameter
$b$, and sometimes we use the notation $\G_b$.
We call any tree $\G_b$, with an arbitrary $b$,
\textsl{homogeneous}.

The metric topology and the Lebesgue measure on $\G$ are
introduced in a natural way.
The space
$\plainL2(\G)$ is understood as
$\plainL2$ with respect to this
measure.

A function $f$ on $\G$ belongs to the
Sobolev space $\plainH1(\G)$ if and
only if it is continuous,
$f\res e\in \plainH1(e)$ for each edge $e$, and
\begin{equation*}
\Vert f\Vert^2_{\plainH1(\G)}
:=
\int_\G(|f'|^2+|f|^2)dx<\infty.
\end{equation*}
As usual, $\plainH{1,0}(\G) = \{f\in \plainH1(\G):f(o)=0\}$.

We define the Dirichlet
Laplacian $-\boldsymbol{\D}$ on $\G$ as the self-adjoint
operator in
$\plainL2(\G)$, associated with the
quadratic form $\int_\G|f'|^2dx$ considered
on the form domain $\plainH{1,0}(\G)$. It is easy to
describe the operator domain
$\dom(\boldsymbol{\D})$ and the action of $\boldsymbol{\D}$.
Evidently $f\in\dom(\boldsymbol{\D})\Rightarrow f\res e\in \plainH2(e)$
for each edge $e$ and the Euler-Lagrange
equation reduces on $e$ to
$\boldsymbol{\D} f=f''$.
In order to describe the matching conditions
at a vertex $v\not=o$,
denote by $f_-$ the restriction $f\res e_v^-$
and by $f_j,\;j=1,\ldots,b$
the restrictions $f\res e_v^j$.
The matching conditions at $v$ are
\begin{equation*}
f_-(v)=f_1(v)=\ldots=f_b(v);\ \
f_1'(v)+\ldots+f_b'(v)=f'_-(v)
\end{equation*}
where the derivatives on each edge are taken in the direction
consistent with the ordering on $\G$.
The first matching condition comes
from the requirement $f\in \plainH1(\G)$
which includes the continuity of $f$,
and the second appears as the
natural condition in the sense of
Calculus of Variations.
At the root $o$ we have the boundary
condition $f(o)=0$. It is easy to check that the conditions listed are
also sufficient for $f\in\dom (\boldsymbol{\D})$.

Along with the Laplacian $-\boldsymbol{\D}$ we
shall be interested also in the
Schr\"odinger operators with a
real, bounded and \textsl{symmetric}
(that is, depending
only on $|x|$) potential $V$:
\begin{equation}\label{lap:1}
\BA_Vf:=-\boldsymbol{\D} f+V(|x|)f,\qquad f\in\dom(\boldsymbol{\D}).
\end{equation}
The operator $\BA_V$ is self-adjoint.
Its quadratic form is given by
\begin{equation*}
\ba_V[f]=\int_\G(|f'|^2+V(|x|)|f|^2)dx,\ \
f\in \plainH{1,0}(\G).
\end{equation*}

\subsection{ The orthogonal decomposition of $\plainL2(\G)$}
Our techniques is based upon the orthogonal decomposition
of $\plainL2(\G)$ into a family of
subspaces associated with a class of
subtrees of $\G$.
Given a subtree $T\subset\G$, we say that a function
$f\in \plainL2(\G)$ belongs to
$\CF_T$ if and only if
\begin{equation*}
f=0\qquad\text{outside}\;T
\end{equation*}
and
\begin{equation}\label{lap:2}
f(x)=f(y)\qquad\text{if}\ x,y\in T\ \text{and}\ |x|=|y|.
\end{equation}
Evidently $\CF_T$ is a closed subspace
of $\plainL2(\G)$. It is easy to
describe the operator $\CP_T$ of
orthoprojection onto $\CF_T$. To this end,
introduce the function
\begin{equation*}
b_T(t)=\#\{x\in T:|x|=t\}.
\end{equation*}
In particular,
\begin{equation}\label{lap:2c}
b_\G(t)=b^k\ \ \text{for}\;k-1\le t<k,\;k\in\mathbb N.
\end{equation}
It is clear that
\begin{equation*}
(\CP_Tf)(x)=
\begin{cases}
(b_\G(|x|))^{-1}
\sum\limits_{y\in T:\,|y|=|x|}  f(y)& \text{for}\; x\in T;\\
0 & \text{for}\;x\notin T.
\end{cases}
\end{equation*}

We shall need the subspaces $\CF_T$,
associated with the subtrees of
two following types. Given a vertex $v$, let
\begin{equation*}
T_v=\{x\in\G:x \succeq v\}.
\end{equation*}
Given an edge $e=\langle v,w\rangle$, let
\begin{equation*}
T_e=e\cup T_w.
\end{equation*}
In particular, $T_{e_{_0}}=T_o=\G$.
For the sake of brevity, for any $v\ne o$
below we use the notation
$\CF_v,\;\CF_v^j$ for $\CF_{T_v},\;
\CF_{T_{e_v^j}}$.
It is clear that the subspaces
$\CF_v^1,\ldots,\CF_v^b$ are mutually
orthogonal and their orthogonal sum
\begin{equation*}
\widetilde{\CF_v}=\CF_v^1\oplus\ldots\oplus\CF_v^b
\end{equation*}
contains $\CF_v$. Denote by
$\CF'_v = \widetilde{\CF_v}\ominus\CF_v$ the orthogonal
complement.
The next theorem is a direct
consequence of \cite{NS}, Theorem 5.1 and
Lemma 5.2,
where a more general class of trees was considered.
Later the result was re-discovered by
R.Carlson \cite{carlson2}, in a slightly different setting.
A new detailed exposition,
most convenient for our
purposes, was recently given in \cite{Sol}.

\begin{thm}\label{lap:thm1} Let $\G=\G_b$ for some $b>1$.
\begin{itemize}
\item[(i)]
The subspaces $\CF'_v$, $o\not=v\in\CV(\G)$ are
mutually orthogonal and orthogonal to $\CF_\G$.
Moreover,
\begin{equation}\label{lap:10}
\plainL2(\G)=\CF_\G\oplus
\sum_{v\in\CV(\G)\setminus \{o\}}\oplus\CF'_v.
\end{equation}
\item[(ii)]
Let $V(t)$ be a real, measurable and bounded function on $\R_+$.
Then the decomposition \eqref{lap:10} reduces the Schr\"odinger
operator \eqref{lap:1}, and in particular
the Laplacian $-\boldsymbol{\D}=\BA_0$.
\end{itemize}
\end{thm}

\bigskip

\subsection{Parts of $\BA_V$ in
the subspaces $\CF_\G$, $\CF'_v$}
According to \thmref{lap:thm1}, the description of the spectrum
$\s(\BA_V)$
reduces to the similar
problem for the parts of $\BA_V$ in the
components of the decomposition
\eqref{lap:10}. Consider at first the part
of $\BA_V$
in the subspace $\CF_\G$.
It is more convenient (and equivalent)
to deal with the quadratic form $\ba_V$.

It is natural to identify a function
$f\in \CF_\G$ with the function
$\varphi$ on $\R_+$, such that $\varphi(t)=f(x)$ for $|x|=t$.
The operator $\Pi: f\mapsto\varphi$
acts as an isometry of $\CF_\G$
onto the weighted space $\plainL2(\R_+,b_\G)$
with the norm given by
\begin{equation*}
\Vert\varphi\Vert^2_{\plainL2(\R_+,b_\G)}
=\int_{\R_+}|\varphi(t)|^2b_\G(t)dt.
\end{equation*}
Then
\begin{equation}\label{lap:12}
\ba_V[f] =
\int_{\R_+}\bigl(|\varphi'(t)|^2
+ V(t)|\varphi(t)|^2\bigr)b_\G(t)dt,\ \varphi = \Pi f.
\end{equation}
Its domain is the weighted Sobolev
space $\plainH{1,0}(\R_+,b_\G)$ whose
norm is defined by the quadratic
form \eqref{lap:12} with $V\equiv1$.
The corresponding operator
$\BA_V\res\CF_\G$ turns into an
operator acting in $\plainL2(\R_+,b_\G)$.
It is not difficult to describe it
explicitly, however it is more
natural to pass on to the
operators acting in the ``usual''
$\plainL2(\R_+)$. To this end we make the
substitution
\begin{equation}\label{lap:13}
y(t)=b_\G(t)^{1/2}\varphi(t).
 \end{equation}
Then
\begin{equation*}
\Vert y\Vert^2_{\plainL2(\R_+)}
=\Vert\varphi\Vert^2_{\plainL2(\R_+,b_\G)}.
\end{equation*}
Since $b_\G(t)$ is a step function, we also have
\begin{equation*}
a_V[y] := \ba_V[f] =
\int_{\R_+}\bigl(|y'(t)|^2+V(t)|y(t)|^2\bigr)dt.
\end{equation*}
However, the domain of $a_V$ does not
coincide with $\plainH{1,0}(\R_+)$,
since the function $y(t)$ may have
jumps at the points $n\in\mathbb N$.
More exactly, it follows from
\eqref{lap:2c} and \eqref{lap:13} that
$\dom(a_V)$ consists of
functions
\begin{equation}\label{lap:n}
y\in \plainH1(0,1)\times \plainH1(1,2)\times
\ldots\times \plainH1(n-1,n)\times\ldots
\end{equation}
such that
\begin{equation}\label{lap:16}
y(0)=0;\qquad y(n+)=b^{1/2}y(n-), \   \forall n\in\mathbb N,
\end{equation}
and
\begin{equation}\label{lap:17}
\int_{\R_+}\bigl(|y'(t)|^2+|y(t)|^2\bigr)dt<\infty.
\end{equation}
The self-adjoint operator in $\plainL2(\R_+)$, associated
with this quadratic form,
on each interval $(n-1,n)$, $n\in\mathbb N$ acts as
\begin{equation*}
A_V y=-y''+V(t)y.
\end{equation*}
Its domain $\dom( A_V)$ consists of all functions
\begin{equation*}
y\in \plainH2(0,1)\times \plainH2(1,2)\times
\ldots\times \plainH2(n-1,n)\times\ldots
\end{equation*}
satisfying the conditions \eqref{lap:16} and
\begin{equation}\label{lap:19}
y'(n+)=b^{-1/2}y'(n-),\  \forall n\in\mathbb N,
\end{equation}
and also
\begin{equation*}
\int_{\R_+}\bigl(|y''(t)|^2+|y(t)|^2\bigr)dt<\infty.
\end{equation*}
(Here and in \eqref{lap:17} it
would be more accurate to write
$\sum\limits_{n=1}^\infty\int\limits_{n-1}^n $ rather than
$\int\limits_{\R_+}$.)
So we have proved the following

\begin{lem}\label{lap:lem2}
The part of the operator $\BA_V$ in the subspace $\CF_\G$
is unitarily equivalent to the operator $A_V$ in $\plainL2(\R_+)$.
\end{lem}

\bigskip

Now we turn to the operators
$\BA_V\res\CF'_v$, $v\neq o$.
It follows
from the symmetry properties of
the tree $\G$ and of the potential $V(|x|)$
that all such operators with the same value
of $\gen(v)=k$ can be identified
with each other.
In order to reduce them to the operators in
$\plainL2(\R_+)$, introduce the
``shifted'' potentials
\begin{equation}\label{lap:trun}
V_k(t)=V(t+k),\ \ t>0,\ \ k=0,1,\ldots.
\end{equation}
In particular, $V_0=V$.

\begin{lem}\label{lap:lem3}
Let $\G=\G_b$ and $v\in\CV(\G),\ \gen(v)=k>0$.
Then
the operator $\CA_V\res\CF'_v$ is unitarily
equivalent to the orthogonal
sum of $(b-1)$ copies of the operator
$A_{V_k}$.
\end{lem}

For the formal proof,
see \cite{Sol}.
On the qualitative level,
the result follows from
the fact that the restriction of the operator $\BA_V$
to the subspace $\widetilde{\CF_v}$
reduces to orthogonal sum of $b$ copies
of the operator $A_{V_k}$. The passage
to the subspace $\CF'_v$ corresponds to
the withdrawal of one of these copies.

\subsection{ The orthogonal
decomposition of the operators $\BA_V$}
Now we are in position to present
the final result of this section.
Below
$A^{[r]}$ stands for the orthogonal
sum of $r$ copies of a self-adjoint
operator $A$.

\begin{thm}\label{lap:th4} Let $\G=\G_b$
and the function $V$ be
real, measurable and bounded on $\R_+$. Then
the Schr\"odinger operator \eqref{lap:1}
on $\G$ is unitarily equivalent
to the orthogonal sum of the operators acting in
$\plainL2(\R_+)$:
\begin{equation*}
\BA_V\sim A_V\oplus\sum_{k\in\mathbb N}
\oplus {A_{V_k}}^{[b^{k-1}(b-1)]}.
\end{equation*}
In particular, for the Laplacian $-\boldsymbol{\D}=\BA_0$ we get
\begin{equation*}
-\boldsymbol{\D}\sim {A_0}^{[\infty]}.
\end{equation*}
\end{thm}

This Theorem is a direct consequence of
Theorem \ref{lap:thm1} and
Lemmas \ref{lap:lem2},  \ref{lap:lem3}, if one remembers
that the total number of vertices of generation $k$
equals $b^{k-1}$.

\section{Spectrum of the Laplacian on $\G_b$}
\label{tree3:sect}

\subsection{The operator $A$ on the whole line}
Along with the operator $A_0$
in $\plainL2(\R_+)$ defined as $A_0y=-y''$
with the boundary and matching conditions \eqref{lap:16} and
\eqref{lap:19}, consider the
similar operator, say $A$, in $\plainL2(\R)$:
\begin{equation*}
(Ay)(t)=-y''(t),\ \ t\not\in\Z,
\end{equation*}
on the analogous domain supplied with the matching conditions
\begin{equation}\label{sp:1}
y(n+) = b^{1/2} y(n-),\ \ y'(n+) = b^{-1/2}y'(n-),\ \
n\in\Z.
\end{equation}
The spectrum of $A$ can be found
by means of the standard Floquet procedure.
The related quasi-periodic problem is
\begin{equation*}
y'' + \mu^2 y = 0,\ \ y(1+) = e^{i\xi} y(0+),\ \
y'(1+) = e^{i\xi} y'(0+),
\end{equation*}
with the parameter (quasi-momentum) $\xi\in [0, 2\pi)$.
Taking into account the matching conditions at the point $n=1$,
we can re-write this as
\begin{gather}
y''(t) + \mu^2 y(t) = 0, \ 0<t<1;\label{sp:eq}\\
y(1-) = b^{-1/2} e^{i\xi} y(0+),\ \
y'(1-) = b^{1/2} e^{i\xi} y'(0+).\notag
\end{gather}
It is quite straightforward to
calculate the eigenvalues of the
 problem \eqref{sp:eq}.
Introduce the function
\begin{equation}\label{sp:phi}
\varphi(\xi) = \arccos \ \frac{\cos\xi}{R},\ \
R = \frac{b^{1/2}+ b^{-1/2}}{2} > 1.
\end{equation}
Then the numbers $\mu_l$
(square roots of eigenvalues) are given by
\begin{equation}\label{sp:mu}
\mu_l(\xi) =
\begin{cases}
\pi(l-1) + \varphi(\xi),\ l \ \textup{is odd},\\
\pi l - \varphi(\xi),\ \ l \ \textup{is even},
\end{cases}
l\in \mathbb N.
\end{equation}
The function $\varphi$ is one-to-one on
the interval $[0, \pi]$. Later we shall also need its inverse:
\begin{equation}\label{sp:6}
\psi(\mu) = \arccos(R\cos\mu),\qquad
\mu\in [\varphi(0), \ \varphi(\pi)] = [\t, \pi - \t].
\end{equation}
where
\begin{equation*}
\t = \arccos(1/R).
\end{equation*}
It follows easily from \eqref{sp:phi} that
\begin{equation}\label{effmass:eq}
\begin{cases}
\psi(\mu) = &\ 2^{1/2} (R^2-1)^{1/4} (\mu-\t)^{1/2}
+ O(\mu-\t),\ \mu\to \t+,\\[0.2cm]
\psi(\mu) = &\ \pi - 2^{1/2} (R^2-1)^{1/4} (\mu-\pi+\t)^{1/2}\\
&\  + O(\mu-\pi+\t),\ \mu\to \pi-\t-.
\end{cases}
\end{equation}
Define the segments (``bands'')
\begin{equation*}
{\gb}_l =\bigcup_\xi \mu^2_l(\xi)=
\bigl[(\pi(l-1)+\t)^2,(\pi l-\t)^2\bigr],\ \ l\in \mathbb N
\end{equation*}
and the intervals (``gaps'')
 \begin{equation}\label{sp:5}
\gl_0 = (-\infty, \t^2),\ \ \
\gl_l = \bigl((\pi l - \t)^2,
(\pi l + \t)^2\bigr),\ \ l\in \mathbb N.
\end{equation}
The gaps are labelled so that
$\gl_l$ separates the bands ${\gb}_l$
and  ${\gb}_{l+1}$. The following statement
is a direct consequence of
the Floquet theory.

\begin{lem}\label{sp:specA}
The spectrum of $A$ coincides with the union
of the bands ${\gb}_l,\ l=1,2,\ldots$.
On this set the spectrum is of
the Lebesgue type and of multiplicity two.
\end{lem}

We shall need also the spectral
decomposition of the operator $A$. To this
end, note that
\begin{gather*}
\zeta_l(t,\xi) = c_l(\xi)
\bigl(\cos(\mu_l(\xi)(1-t))
- b^{1/2}e^{i\xi}\cos(\mu_l(\xi)t)\bigr),
\\
c_l(\xi)=\sqrt{2(b+1)^{-1}}|\sin\mu_l(\xi)|^{-1},\
0 < t <1,
\end{gather*}
is the normalized in $\plainL2(0,1)$ eigenfunction of the equation
\eqref{sp:eq} corresponding to the eigenvalue $\mu^2_l(\xi)$. It
follows from \eqref{sp:mu}, \eqref{sp:phi} that
$\z_l(t, \xi)$
is smooth in $\xi$ on each band $\gb_l$. Let us extend each
function $\zeta_l(t,\xi)$ to all $t\in\R$ in the following way.
Let $\om_l(t,\xi)$ be the periodic (in $t$) extension of the
function $e^{-it\xi}\zeta_l(t,\xi)$ from the interval $[0,1)$ to
$\R$. Then we define $\zeta_l(t,\xi)$ on the whole of $\R$ by the
equation
\begin{equation*}
\zeta_l(t,\xi)=e^{it\xi}\om_l(t,\xi),\ \ t\in\R.
\end{equation*}

Let $P_l$ be the spectral projection of
$A$ associated with the band $\gb_l$.
The map
\begin{equation*}
(U_ly)(\xi)=(2\pi)^{-{1/2}}\int_\R\overline{\zeta_l(t,\xi)}y(t)dt
\end{equation*}
defines the unitary operator from $\plainL2(\R)$ onto
$\plainL2(-\pi,\pi)$
which diagonalizes $AP_l$, namely
\begin{equation*}
(U_lAP_ly)(\xi)=\l_l(\xi)(U_lP_ly)(\xi),\ \ \l_l(\xi)=\mu^2_l(\xi).
\end{equation*}
The adjoint operator
$U^*:\plainL2(-\pi,\pi)\to\plainL2(\R)$ is given by
\begin{equation*}
(U_l^*z)(t) = (2\pi)^{-{1/2}}
\int_{-\pi}^{\pi}\zeta_l(t,\xi)z(\xi)d\xi
=(2\pi)^{-{1/2}}\int_{-\pi}^{\pi}e^{it\xi}
\om_l(t,\xi)z(\xi)d\xi.
\end{equation*}
Denoting by $[m]$ the
operator of multiplication by a scalar function $m$,
we get the spectral decomposition of $A$ in the form
\begin{equation}\label{sp:dec}
A=\sum_{l\in\mathbb N}U_l^*[\l_l]U_lP_l.
\end{equation}

\subsection{Spectrum of the operators $A_0$ and $\BA_0$}

\begin{thm}\label{sp:thmA}
The spectrum of the operator
$A_0$ consists of the bands
${\gb}_l,\ l\in\mathbb N$ and of
the simple eigenvalues
$\l_l=(\pi l)^2, \ l\in \mathbb N$.
The corresponding eigenfunctions
(normalized in $\plainL2(\R_+)$) are
\begin{gather}
y_l(t)=c(b)b^{-{n/2}}\sin (\pi lt),\ \
t\in(n-1,n),\  n\in\mathbb N
\label{sp:n5}\\
c(b)=(2(b-1))^{-1/2}.\notag
\end{gather}
\end{thm}
\begin{proof}
The operator $A_0$ is non-negative, so its spectrum lies
on $[0,\infty)$.

\vskip0.3cm

1. BANDS.
Let $D$ be
the operator in
$\plainL2(\R_+)$, defined as follows:
its operator domain coincides with the
quadratic domain of $A_0$, i.e.
is defined by \eqref{lap:n} -- \eqref{lap:17},
and for $y$ from this domain
\begin{equation}\label{sp:20}
(Dy)(t)=-iy'(t),\ t\notin\mathbb N.
\end{equation}
The operator $D$ is closed and its
adjoint $D^*$ acts by the
same formula \eqref{sp:20} on the domain consisting of those
functions $y$ from
the direct product \eqref{lap:n} which satisfy \eqref{lap:17} and
the matching conditions
similar to the ones in \eqref{lap:16} but with
the factor
$b^{1/2}$ replaced by $b^{-1/2}$; there is no boundary condition at
 $t=0$.

It is easy to see that $A_0 = D^*D$.
Along with $A_0$, consider the operator $DD^*$.
It acts as
$(DD^*y)(t)=-y''(t),\ t\not\in\mathbb N$, and
its domain is described by the boundary condition $y'(0)=0$
and the matching conditions
\begin{equation*}
y(n+) = b^{-1/2} y(n-),\ \ y'(n+) = b^{1/2}y'(n-),\ \
n\in\mathbb N.
\end{equation*}
According to the general operator theory,
the non-zero spectra of the operators
$A_0=D^*D$ and $DD^*$ coincide.

Now, in the definition of the operator
$A$ let us replace the matching
condition at $t=0$ by the boundary conditions
\begin{equation*}
y(0+)=0,\ \ y'(0-)=0.
\end{equation*}
The new operator, say $A'$, splits into the orthogonal sum,
$A'=A_0\oplus A'_0$ where the operator
$A'_0$ acts in $\plainL2(\R_-)$.
Its description is clear from the
construction and it is easy to see
that the substitution $t\mapsto -t$
reduces $A'_0$ to $DD^*$. The essential
spectrum of $A'$ is the same as that of $A$, i.e.
$\cup_{l\in\mathbb N}\gb_l$. It also coincides with the union
of the essential spectra of the operators $A_0$ and $A'_0$, i.e.
with each of them. It follows that the
essential spectrum of $A_0$
coincides with the spectrum
of $A$ given by \lemref{sp:specA}.

\vskip0.3cm

2. EIGENVALUES. The fact that each
function $y_l(t)$, cf. \eqref{sp:n5},
is the eigenfunction corresponding to
the eigenvalue $(\pi l)^2$, can be
verified by the direct inspection.
Any two solutions satisfying the
boundary condition $y(0)=0$ are
proportional to each other, so that
this eigenvalues are simple.
The direct inspection shows also that
$\l=0$ is not an eigenvalue.
So it remains to show that any number
$\l=k^2>0$ with $\pi^{-1}k\not\in\mathbb N$
can not be an eigenvalue.
For this purpose we use the explicit formulae for the solutions
of the equation
\begin{equation}\label{sp:eqno}
y''(t)+k^2 y(t)=0,\ \ t\not\in\mathbb N
\end{equation}
under the matching conditions \eqref{sp:1}.
Namely, let $q_1,q_2$ be found
from the quadratic equation
\begin{equation}\label{sp:n3}
q^2-2Rq\cos k+1=0
\end{equation}
where $R$ is defined in \eqref{sp:phi}. Suppose that $q_1\neq q_2$,
that is $R|\cos k|\neq1$. The functions
\begin{gather}
y_j(t)=\bigl(b^{1/2}\sin k(n-t)+q_j\sin k(t-n+1)\bigr)q_j^{n-1},
\label{sp:sol}\\
 n-1<t<n,\ \ n\in\mathbb N,\ \ j=1,2\notag
\end{gather}
are solutions of the problem \eqref{sp:eqno} -- \eqref{sp:1}.
Their Wronskian is equal to
$y_1y'_2-y'_1y_2=b^{1/2}(q_2-q_1)k\sin k$,
so that the solutions $y_1,\ y_2$
are linearly dependent only if
$\pi^{-1}k\in\mathbb N$ which is the excluded case.
Any solution satisfying
the condition $y(0+)=0$ is proportional to the function
\begin{gather}
y_0(t)=\frac{y_2(t)-y_1(t)}{q_2-q_1}\notag\\
=
b^{1/2}\frac{q_2^{n-1}-q_1^{n-1}}{q_2-q_1}\sin k(n-t)+
\frac{q_2^n-q_1^n}{q_2-q_1}\sin k(t-n+1),\label{sp:sol1}\\
 n-1<t<n,\ \ n\in\mathbb N.\notag
\end{gather}
For $\pi^{-1}k\not\in\mathbb N$ this function does not lie in
$\plainL2(\R_+)$ and hence, is not an eigenfunction.
If $R|\cos k|=1$, then $q_1=q_2=\pm1$
and it is easy to see that
there also are no $\plainL2$-solutions of the problem
\eqref{sp:eqno} -- \eqref{sp:1}, and we are done.
\end{proof}

The result for the operator $\BA_0$,
that is for the Laplacian on the tree,
immediately follows from \thmref{lap:th4} and  \thmref{sp:thmA}.

\begin{thm}\label{sp:thmA1}
The spectrum of the operator $\BA_0$ is of infinite multiplicity
and consists of the bands ${\gb}_l$ and the eigenvalues
$\l_l=(\pi l)^2, \ l\in \mathbb N$.
\end{thm}

We see that the gap $\gl_0$ of the operator
$A$ is also the gap for $\BA_0$,
and each gap $\gl_l$ of $A$ with $l\ge1$
splits into two gaps when we
turn to the operator $\BA_0$:
\begin{equation*}
\gl_{l,-}=\bigl((\pi l-\t)^2,(\pi l)^2\bigr);\ \  \gl_{l,+}=
\bigl((\pi l)^2,(\pi l+\t)^2\bigr).
\end{equation*}

\subsection{Global quasi-momentum and density of states}
Define the density of states for
the operators $A$ and $A_0$ as the limit
\begin{equation}\label{sp:9}
\rho(\l) = \lim \frac{N_P(\l; \Delta)}{|\Delta|},\ \
|\Delta|\to\infty.
\end{equation}
Here we denote $\Delta = (0, L),\ L\in \mathbb N$,
and $N_P(\l) = \#\{j: \mu_l^2 < \l \}$
is the counting function for
the operator $By=-y''$ which at the points $1,\dots,L-1$
has the same matching conditions
as in \eqref{sp:1}, and also satisfies the boundary conditions
\begin{equation*}
y(0) = b^{1/2} y(L),\ \
y'(0) = b^{-1/2} y'(L).
\end{equation*}
The subscript $P$ in the notation for the counting function
indicates that the operator
$B$ has the boundary conditions
of this type.
If the limit \eqref{sp:9} exists for these conditions, then
it will also exist for any other conditions, and its value
will not depend on them.
Later, in order to calculate
the density of states we shall use the same
formula \eqref{sp:9}, but with the counting function of
the Dirichlet problem. In this case we do not
use any subscripts and simply write $N(\l)$.

Let us find eigenvalues of $B$.
Denote $k = \sqrt\l$, then choose solutions on every interval
$(n, n+1)$ in the form
\begin{equation*}
y(t)=\a_n\cos k(t-n)+\b_n\sin k(t-n).
\end{equation*}
In view of the matching conditions,
we come, with the notations $c = \cos k, s = \sin k$,
to the equalities
\begin{equation}\label{sp:n2}
\a_n=b^{1/2}(\a_{n-1}c+\b_{n-1}s),\ \
\b_n=b^{-{1/2}}(-\a_{n-1}s+\b_{n-1}c)
\end{equation}
for $n = 0, \dots, L$.
Here we have identified the points with $n=0$ and $n = L$,
so that $\a_0 = \a_L$ and $\b_0 = \b_L$.
To solve this system
introduce the functions
\begin{equation*}
\CA(z) = \sum_{n=0}^{L-1} \a_n z^n,\ \
\CB(z) = \sum_{n=0}^{L-1} \b_n z^n,
\end{equation*}
where $z$ runs over the set of all complex
numbers such that $z^L = 1$. Then by \eqref{sp:n2}
\begin{equation*}
\CA(z) = b^{1/2} z \bigl(c\CA(z) + s \CB(z)\bigr), \ \
\CB(z) = b^{-1/2} z \bigl(-s\CA(z) + c \CB(z)\bigr).
\end{equation*}
This system of two equations has non-trivial solution iff
its determinant is identically zero:
\begin{gather*}
\det
\begin{pmatrix}
cb^{1/2}z-1 & sb^{1/2}z\\
-sb^{-1/2}z & cb^{-1/2}z - 1
\end{pmatrix}\\
= c^2 z^2 - cb^{-1/2}z - cb^{1/2}z + 1 + s^2 z^2
= z^2 - 2cR z + 1 = 0,
\end{gather*}
whence
\begin{equation*}
R\cos k =Rc= \frac{z+z^{-1}}{2} = \cos\frac{2\pi n}{L},
\ \ n = 0, 1, \dots, L-1.
\end{equation*}
It is convenient to  write the formulae for the eigenvalues
in terms
of the function
$\varphi$ defined by \eqref{sp:phi}, and the formulae for
$N_P(\l)$ -- in terms of the ``global quasi-momentum''
$\om(\l)$ which
we now define.
Namely, $\om(\l)=\pi l$ if $\l\in\gl_l$, and
for $\l\in \gb_l$
\begin{equation}\label{sp:22}
\om(\l) =
\begin{cases}
\pi (l-1) + \psi\bigl(\sqrt\l-\pi(l-1)\bigr),\ & l \ \
\textup{is odd},\\
\pi l - \psi\bigl(\pi l - \sqrt\l \bigr),\ & l \ \ \textup{is even}.
\end{cases}
\end{equation}
Here $\psi$ is the function inverse to
$\varphi$, cf. \eqref{sp:6}.
Evidently
\begin{equation}\label{sp:7a}
c\sqrt{\l-\t^2}\le \om(\l)\le C\sqrt\l,\ \ \l\ge0.
\end{equation}
The eigenvalues of the operator $B$ are given by the formulae
\begin{gather*}
\mu_l(j) =\pi(l-1) + \varphi\biggl(\dfrac{2\pi j}{L}\biggr),
\ l\ \textup{odd};\ \
\mu_l(j) = \pi l -   \varphi\biggl(\dfrac{2\pi j}{L}
\biggr),\ l\ \textup{even},\\
l\in \mathbb N,\ j = 0, 1, \dots, L-1.
\end{gather*}
The number $N_P(\l)$ depends on the location of $\l$.
For instance, if $\l\in \gl_l$, then
$N(\l) = lL$, so that
\begin{equation*}
\rho(\l) = \frac{N_P(\l, \Delta)}{|\Delta|}
= l=\frac{1}{\pi}\om(\l),\ \ \l\in \gl_l.
\end{equation*}
To cover the case $\l\in \gb_l$ we shall consider
two options: $l$ is odd or $l$ is even.
Suppose first that $l$ is odd and denote
$\sqrt\l = \pi(l-1) + \chi$. Then
\begin{gather*}
\frac{N_P(\l; \Delta)}{L}
=  l-1 + \ell_1(\chi),\\[0.2cm]
\ell_1(\chi) = \frac{1}{L}\#\left\{
j \in [0, L-1): \varphi\biggl(\dfrac{2\pi j}{L}
\biggr)< \chi
\right\}.
\end{gather*}
The term $\ell_1(\chi)$ can be estimated as follows:
\begin{equation*}
\left|\ell_1(\chi) - \frac{1}{\pi}\psi(\chi)\right|\le 2L^{-1}.
\end{equation*}
Hence by \eqref{sp:22}
\begin{equation}\label{sp:11}
\biggl|\frac{N_P(\l; \Delta)}{L}
- \frac{1}{\pi}\ \om(\l)
\biggr|\le 2L^{-1},\ \ \forall \l > 0.
\end{equation}

Suppose now that $l$ is even.
Denote $\sqrt\l  = \pi l - \chi$. Then
\begin{gather*}
\frac{N_P(\l; \Delta)}{L} =
l - \ell_2(\chi),\\[0.2cm]
\ell_2(\chi) = \frac{1}{L}\#\left\{ j \in [0, L-1):
\varphi\biggl(\dfrac{2\pi j}{L} \biggr)> \chi \right\}.
\end{gather*}
The term $\ell_2(\chi)$ can be estimated as follows:
\begin{equation*}
\left|\ell_2(\chi) + \frac{1}{\pi}\psi(\chi) - 1\right|\le 2L^{-1}.
\end{equation*}
Using \eqref{sp:22} again, we get \eqref{sp:11}. All this results in the
formula
\begin{equation}\label{sp:11a}
\rho(\l)=\frac{1}{\pi}\,\om(\l)
\end{equation}
which is well known for the clasical Hill operator.
\vskip0.3cm
Relying upon the estimate \eqref{sp:11} we shall prove a
similar estimate for the counting function
of the Dirichlet problem on an arbitrary interval
$(R_1,R_2)$, not necessarily with integer $R_1,,R_2$.

\begin{thm}\label{density:thm}
Let $\Delta =(R_1,R_2),\ \ R_1,R_2\in\R,\ \  R_1<R_2$.  Then
the inequality holds:
\begin{equation}\label{sp:12}
\left| \frac{N(\l; \Delta)}{|\Delta|} - \rho(\l)\right|
\le C\frac{1+\sqrt\l}{|\Delta|}
\end{equation}
for all $\l>0$,
with a universal constant $C$.
\end{thm}

\begin{proof}
It is well known that for the Dirichlet realization of the operator
$-y''$ on $\Delta$ (with no matching conditions inside!) the counting
function is controlled by $C|\Delta|\sqrt\l$, with a universal constant.
Since the number of integer points  inside $\Delta$ is
not greater than
$|\Delta|+1$, it follows from
the decoupling principle that
\begin{equation}\label{sp:13}
{N(\l; \Delta)}\le C(|\Delta|\sqrt\l +|\Delta|+1).
\end{equation}
If the length of $\Delta$ is small, say $|\Delta|\le2$, then
\eqref{sp:12} is implied by \eqref{sp:7a} and \eqref{sp:13}.

Let now $|\Delta|>2$.
Without loss of generality, we may assume $0\le R_1<1$.
Define $L=[R_2]$ (the integer part of $R_2$), then $L\ge2$.
Let $\Delta_+$ and $\Delta_-$
be the intervals $(0, L+1)$ and $(1, L)$ respectively.
Then, clearly,
\begin{equation*}
N(\l; \Delta_-)\le N(\l; \Delta)\le N(\l; \Delta_+),
\end{equation*}
by variation argument. Furthermore,
by the decoupling principle,
\begin{equation*}
\bigl|N(\l; \Delta_\pm)
- N_P(\l; \Delta_\pm)\bigr|\le 4,
\end{equation*}
so that
\begin{equation*}
\frac{N_P(\l; \Delta_-)}{|\Delta|} - \frac{4}{|\Delta|}
\le \frac{N(\l; \Delta)}{|\Delta|}
\le \frac{N_P(\l; \Delta_+)}{|\Delta|} + \frac{4}{|\Delta|}.
\end{equation*}
Therefore
\begin{equation*}
\biggl|\frac{N(\l; \Delta)}{|\Delta|} - \rho(\l)\biggr|
\le \max_{\pm}\biggl|\frac{N_P(\l; \Delta_\pm)}{|\Delta|} -
\rho(\l)\biggr|
+ \frac{4}{|\Delta|}.
\end{equation*}
Let us estimate the r.h.s. with the "$-$" sign. The modulus
equals
\begin{equation*}
\biggl|\frac{N_P(\l; \Delta_-)}{L-1}\frac{L-1}{|\Delta|}
- \rho(\l)\biggr|
\le \biggl|\frac{N_P(\l; \Delta_-)}{L-1}
- \rho(\l)\biggr| + \frac{2 N_P(\l; \Delta_-)}{|\Delta|(L-1)}.
\end{equation*}
In view of \eqref{sp:11},
the first term in the r.h.s. is bounded by
$2(L-1)^{-1}\le3L^{-1}$ and in view of
\eqref{sp:13}, the second term
is bounded by
\begin{equation*}
\frac{C((L-1)\sqrt\l+L)}
{|\Delta|(L-1)}\le\frac{2C(\sqrt\l+1)}{|\Delta|}.
\end{equation*}
Repeating the same argument for the "$+$" sign,
we arrive at \eqref{sp:12}.
\end{proof}

\vskip0.3cm

We conclude this section by
discussing the H\"older properties of the
global quasimomentum $\om(\l)$
and thus, those of the density of states
$\rho(\l)$.
It is clear from \eqref{effmass:eq}
that near the edges of the gap $\gl_l = (\l_-, \l_+)$
the function $\rho$ has the following behaviour:
\begin{equation*}
\rho(\l) = \rho(\l_{\pm})\pm \frac{\sqrt2}{\pi}
\left[
\frac{(R^2-1)}{\l_{\pm}}
\right]^{\frac{1}{4}}
(\l-\l_{\pm})^{\frac{1}{2}}
+ O(\l-\l_{\pm}),\  \l\to \l_{\pm}+ 0\pm.
\end{equation*}
Together with the formula \eqref{sp:7a} this asymptotics
guarantees that
\begin{equation}\label{effmass1:eq}
|\rho(\l) - \rho(\l_{\pm})|\ge c|\l-\l_{\pm}|^{\frac{1}{2}},\
\l\in\R,
\end{equation}
with a constant $c$ depending on $l$.
The formula \eqref{sp:phi} also ensures that the function $\psi$ is
$1/2$-H\"older continuous, i.e.
\begin{equation*}
|\psi(\mu_2)-\psi(\mu_1)|\le C|\mu_2-\mu_1|^{1/2},\ \
\mu_1,\mu_2\in [\t,\pi-\t].
\end{equation*}
Using \eqref{sp:7a}, one can immediately
extend this information to the function $\om$:
\begin{equation*}
|\om(\l_2)-\om(\l_1)|\le C\bigl(|\l^{1/2}_2-\l^{1/2}_1|^{1/2}+
|\l^{1/2}_2-\l^{1/2}_1|\bigr), \ \
\l_1,\l_2\ge\t,
\end{equation*}
with a constant $C$ independent of $\l_1, \l_2$.
Later we shall use a less precise, but somewhat more compact
consequence of this estimate and \eqref{sp:11a}:
\begin{equation}\label{holder1:eq}
|\rho(\l_2) - \rho(\l_1)|\le C|\l_2-\l_1|^{1/2},\ \ \l_1, \l_2\in\R,
\end{equation}
with a universal constant $C$.

\section{Operator $\BA_V$
with a decaying potential.
Eigenvalues in the gaps}
\label{tree4:sect}

\subsection{Functions $M(\l)$ and $N(\l_1,\l_2)$}
\label{tree4.1:subsect}
Here we turn to the study of the
spectrum of the Schr\"odinger operators
$\BA_V$, cf. \eqref{lap:1},
with the real-valued and bounded potential $V(|x|)$
which in an appropriate
sense decays as $|x|\to\infty$. The essential
spectrum of $\BA_V$ is the
same as for the unperturbed operator $\BA_0$
(i.e. Laplacian) and therefore,
is given by \thmref{sp:thmA1}. The
spectrum of $\BA_V$ may include
also eigenvalues lying in the gaps of
$\BA_0$. For their study,
the following quantities are standardly used.

Let $\CC$ be a self-adjoint operator
in a Hilbert space, and let $V$ be
its relatively compact perturbation;
we denote $\CC_V=\CC+V$. Suppose that
the interval $(\l_-,\l_+)$
is a gap in $\s(\CC)$.
Let $\l\in(\l_-,\l_+)$. Define the counting
function $M(\l; \CC_V)$
as the number of eigenvalues of $\CC_{\a V}$ crossing
the point $\l$ while $\a$ varies from $0$ to $1$.
In other words,
\begin{equation*}
M(\l;\CC_V)=\sum_{0<\a<1}\dim\ker(\CC+\a V-\l).
\end{equation*}
If $\l$ coincides with one of
the ends of a gap, the function $M(\l; \CC_V)$
is defined as the corresponding
one-sided limit. If $V$ is a perturbation
of fixed sign, that is if $V = \pm q$ with a $q \ge 0$, then
the function $M(\l)$ is increasing (for $V = -q$)
or decreasing (for $V = q$) in $\l\in (\l_-, \l_+)$
and increasing in $q$.

For any subinterval
$(\l_1,\l_2)\subset(\l_-,\l_+)$ the function
$N(\l_1,\l_2; \CC_V)$ is defined
as the total multiplicity of eigenvalues
of the operator $\CC_V$, lying in $(\l_1,\l_2)$. Note that if
$\l_-=-\infty$ and $\l\le\l_+$, then
\begin{equation*}
N(\l; \CC_V):=N(-\infty,\l; \CC_V) = M(\l; \CC_V).
\end{equation*}
According to \thmref{sp:thmA1},
the following equalities hold:
\begin{equation}\label{w:2}
M(\l; \BA_V) = M(\l; A_V)
+ (1-b^{-1})\sum_{k\in\mathbb N}b^kM(\l; A_{V_k}),
\end{equation}
\begin{equation}\label{w:3}
N(\l_1,\l_2; \BA_V)
= N(\l_1,\l_2; A_V)+
(1-b^{-1})\sum_{k\in\mathbb N}b^k N(\l_1,\l_2; A_{V_k}).
\end{equation}
Recall that the potentials
$V_k$ appearing in \eqref{w:2},
\eqref{w:3}
were defined in \eqref{lap:trun}.
These formulae show that the
key step to understanding the behaviour of
the functions $M(\l; \BA_V)$,
$N(\l_1,\l_2; \BA_V)$ consists in studying
the individual terms of the series
\eqref{w:2}, \eqref{w:3}. More precisely,
we need the detailed information
about their behaviour depending on the
parameter $k$.

The study of the sums \eqref{w:2},
\eqref{w:3} is hampered by the presence
of the exponential factors in their r.h.s.
These factors reflect the geometry
of the tree rather than the properties
of the potential $V(t)$. For this
reason, it makes sense to investigate,
along with the functions
$M(\l; \BA_V)$, $N(\l_1,\l_2; \BA_V)$, also the functions
\begin{equation}\label{w:4}
\widetilde M(\l; \BA_V)=\sum_{k\ge0}M(\l; A_{V_k}),
\end{equation}
\begin{equation}\label{w:5}
\widetilde N(\l_1,\l_2; \BA_V)=
\sum_{k\ge0}N(\l_1,\l_2; A_{V_k}).
\end{equation}
For technical reasons, we shall
need also the functions $M$, $N$ for the
operators on intervals $\D\subseteq\R$.
Define $A_{V,\D}$ as the operator in
$L_2(\D)$ acting as
$(A_{V,\D}y)(t)=-y''(t)+V(t)y(t)$ for
$t\not\in\Z$, under the
zero boundary conditions at each
finite end of $\D$ and the matching
conditions \eqref{sp:1}
at the points $n\in\Z\cap\D$. In particular,
$A_{V,\R_+}=A_V$. Often we use
abbreviated notation for the
corresponding
functions $M$, $N$, such as $M(\l,V; \D)$
or even $M(\l;\D)$
when the potential $V$ is fixed.
Note a convenient relation
\begin{equation}\label{shift:eq}
M(\l; A_{V_k}, (R_1, R_2)) = M(\l; A_V, (R_1+k, R_2+k)),
\end{equation}
which is valid for any
$0\le R_1 < R_2\le \infty$ and integer $k$'s.
This formula is useful when it is more natural to
study the dependence of $M$ on the
interval $\D$ than on the potential.

\vskip0.3cm

If $V$ is a function of constant sign,
then there is a useful relationship
between $M(\l; A_{V,\D})$ and the spectrum of the compact operator
\begin{equation}\label{w:6}
T(\l) = T(\l,V,\D)
= |V|^{1/2}(A_{0,\D}-\l I)^{-1}|V|^{1/2},\ \
\l\not\in\s(A_{0,\D}).
\end{equation}
Namely, if $\l$ is a regular point of $A_{0,\D}$, then
\begin{equation}\label{w:7}
M(\l; V,\D)=n_+(1, T(\l,V,\D)),\ \ V\le0;
\end{equation}
\begin{equation}\label{w:8}
M(\l; V,\D)=n_-(1,T(\l,V,\D)),\ \ V\ge0.
\end{equation}
Here $n_\pm(\cdot,T)$ stands for
the counting functions of the positive
and negative eigenvalues
$\pm\l_j^\pm(T)$ of a compact, self-adjoint
operator $T$, that is
\begin{equation*}
n_\pm(s,T)=\#\{j:\l_j^\pm(T)>s\},\ \ s>0.
\end{equation*}
The equalities \eqref{w:7},
\eqref{w:8} proved
very effective in the
problems of the type
considered, see e.g. \cite{Sob}.
Actually, these are facts of
rather general a nature, see e.g. \cite{BirAdv}, Proposition 1.5.

\vskip0.3cm

The following relations have their
prototypes in the theory of the perturbed
Hill operator, see \cite{Sob}, (2.5) -- (2.8). For bounded $\D$
\begin{equation}\label{w:9}
\begin{cases}
M(\l; V,\D) = N(\l; V,\D) - N(\l; 0,\D),\ \ \ V\le 0,\\[0.2cm]
M(\l; V,\D) = N(\l+; 0,\D) -
N(\l+; V,\D),\ \ \ V\ge 0.
\end{cases}
\end{equation}
Further, for any (bounded or unbounded) $\D$
\begin{equation}\label{w:11}
\bigl|N(\l_1,\l_2; V,\D) -
|M(\l_2; V,\D)-M(\l_1; V,\D)|\bigr| 
\le N(\l_1,\l_2; 0,\D)+1.
\end{equation}
One can give a more precise formula: for any
two points $\l_1, \l_2$ such that $N(\l_1, \l_2; 0, \D) = 0$,
we obtain from \eqref{w:9}:
\begin{equation}
\begin{cases}\label{t:eq}
N(\l_1, \l_2; V, \D)
= M(\l_2; V, \D) -  M(\l_1+; V, \D),\\[0.2cm]
N(\l_1, \l_2+; V, \D)
= M(\l_1; V, \D) -  M(\l_2; V, \D).
\end{cases}
\end{equation}
The next two inequalities are usually referred to as the
``decoupling principle''.
Let $\D_1=(a,d)$, $\D_2=(d,c)$,
$-\infty\le a<d<c\le\infty$, and $\D=(a,c)$. Then
\begin{equation}\label{w:12}
\bigl|N(\l_1,\l_2;V,\D)
- \bigl(N(\l_1,\l_2;V,\D_1)
+ N(\l_1,\l_2; V,\D_2)\bigr)
\bigr|\le2,
\end{equation}
\begin{equation}\label{w:13}
\bigl|M(\l; V,\D)-\bigl(M(\l; V,\D_1)+M(\l; V,\D_2)\bigr)
\bigr|\le2.
\end{equation}
The proofs of the relations
\eqref{w:9} -- \eqref{w:13} are either
straightforward, or are
based upon standard facts from the perturbation theory.
Note that the number $1$ rather than $2$ stands in the r.h.s of
the inequalities \cite{Sob}
(2.7) and (2.8) whose analogs are the above
inequalities \eqref{w:12}, \eqref{w:13}.
This difference appears due to the
nature of the matching conditions
at the points $n\in\Z$. If $d\not\in\Z$,
one can replace $2$ by $1$ in \eqref{w:12} and \eqref{w:13}.

\subsection{Individual Weyl asymptotics}
The material presented in this
subsection, is a minor refinement
of \cite{Sob}, Theorems 3.2 and 3.3. We
give it here for the operators we need
in this paper (that is,
the functions in the domains of
the operators considered are subject to
the matching conditions \eqref{sp:1}).
However, it is useful to keep
in mind that the results of
\thmref{w:M}(i) and of \thmref{w:asym}
hold also
for the usual Hill operator.

\vskip0.3cm

For a real-valued
function $V$ on $\R_+$, introduce the quantity
\begin{equation*}
J(V)=
\sum_{n\in\Z} \bigl(\int_{2^{n-1}}^{2^n}
t|V(t)|dt\bigr)^{1/2}.
\end{equation*}
Consider the operator on $\plainL 2(\R_+)$:
\begin{equation}\label{w:15}
K_V y=-y''+Vy,\ \ y\in \plainH2(\R_+),\;y(0)=0.
\end{equation}
In contrast to the operator $A_V$,
the description of $K_V$ involves no
matching conditions, and
the quadratic domain of $K_V$ is
$\plainH{1,0}(\R_+)$.

The following estimate
and asymptotics are particular cases of the results
of \cite{BSolAdv}, Sect. 6;
see also expositions in
\cite{BLap} and \cite{BLapSol}. A close
result was obtained earlier
in \cite{BSolKiev}, Theorems 4.18, 4.19.

\begin{prop}\label{w:est} Let $J(V)<\infty$.
Then the negative spectrum
of the operator $K_V$ is finite and
there exists an absolute constant $C>0$ such that
\begin{equation}\label{w:16}
M(0; K_V)\le C J(V_-).
\end{equation}
Besides, let $g>0$ be the large
parameter and $V_-\not\equiv0$.
The function $M(0; K_{gV})$ satisfies Weyl's asymptotics
\begin{equation}\label{w:17}
\lim g^{-1/2}M(0; K_{gV})= \pi^{-1}
\int_{\R_+}\sqrt{V_-(t)}dt,\ \ g\to\infty.
\end{equation}
\end{prop}
If $|V(t)|$ monotonically decreases, then by
H\"older's inequality
\begin{equation}\label{simplify1:eq}
J(V)/\sqrt6\le \int_{\R_+}|V(t)|^{1/2}dt\le\sqrt6\ J(V).
\end{equation}
Hence, for monotone $|V(t)|$
the function $M(0; K_{gV})$ is controlled
by the r.h.s. of its asymptotics given by \eqref{w:17}.

\vskip0.3cm

Along with $J(V)$, introduce the functional
\begin{equation*}
\widetilde J(V) = \bigl(\int_0^1|V(t)|dt\bigr)^{1/2}+
\sum_{n=1}^\infty \bigl(\int_{2^{n-1}}^{2^n}
t|V(t)|dt\bigr)^{1/2}.
\end{equation*}
Clearly $J(V)\le c\widetilde J(V)$,
therefore in the r.h.s of \eqref{w:16}
$J(V)$ can be replaced by $\widetilde J(V)$.
Compared with \thmref{w:est}, its corollary with
$\widetilde J(V)$
in the r.h.s
ignores the fact that due to
the Dirichlet condition at $0$ the potential
$V$ need not be integrable at this point.
Still, this corollary is
quite
convenient provided one
is dealing with $V$ integrable at $0$.

Present also an estimate for
$M(-1; A_V)$; we need it in the course of
the proof of \thmref{w:M} below.
\begin{equation}\label{w:19}
M(-1; A_V)\le C \T(V):=
C\sum_{n\in\mathbb N}\biggl[\int_{n-1}^n|V(t)|dt\biggr]^{1/2}.
\end{equation}
For the proof, one splits $\R_+$ into the union of the intervals
$(n-1,n]$ and applies to each
interval the well known eigenvalue estimate
for the equation $-y''+y=\l Vy$ with the
Neumann boundary conditions. This
is exactly the way
 in which the
same estimate for $M(-1; K_V)$ was proved in \cite{BirBor}.

It follows from H\"older's inequality that
$\T(V)\le \widetilde J(V)$.
However, the functional $\T(V)$ can not be estimated by $J(V)$.
Note that similarly to \eqref{simplify1:eq},
for a decreasing $|V|$ we have
\begin{equation}\label{simplify:eq}
\widetilde J(V)\le\biggl[\int_0^1|V(t)|dt\biggr]^{1/2}
+ \sqrt{6}\int_{\R_+}
\sqrt{|V(t)|} dt.
\end{equation}

\vskip0.3cm

\begin{thm}\label{w:M}
Let $V$ be a function with a fixed sign and let
$\widetilde J(V)<\infty$.
\begin{itemize}
\item[(i)]
Suppose that $\l\in\overline{\gl_l}$ where
$\gl_l$ is one of the gaps (see \eqref{sp:5}).
Then, given an interval $\D\subseteq\R$,
the estimate
\begin{equation}\label{w:20}
M(\l; V,\D)\le C(\widetilde J(V)+1)
\end{equation}
holds, where the constant $C=C(l)$ does not depend on
$\l\in\bar{\gl_l}$ and $V$.
\item[(ii)]
Suppose in addition that
$\D=\R_+$  (so that $A_{V,\D}=A_V$),
and that
\begin{equation}\label{away:eq}
\l\in I \ \textup{where}\ \
I\subset
\begin{cases}
\overline{\gl_0},\ \ l = 0,\\
\overline{\gl_l}\setminus\{\l_l\},\ l \ge 1,
\end{cases}
\ \ \textup{is a closed interval}.
\end{equation}
Then
\begin{equation}\label{w:201}
M(\l; A_V)\le C'\widetilde J(V)
\end{equation}
with a constant $C'$ uniform in $\l\in I$.
In particular, for $V\le0$
the estimate \eqref{w:201} is uniform in
$\l\le\t^2$.
\end{itemize}
\end{thm}

\begin{proof}
The proof of (i) follows the scheme
suggested in \cite{Sob}.
For this
reason, we only outline
the necessary changes in the argument.
To be definite, we suppose that $V\le0$
and that $l$
(the index of the gap) is even. We start with the
spectral decomposition \eqref{sp:dec}
of the operator $A=A_{0,\R}$
on the whole line.
Set
\begin{gather*}
\l(\xi)=\l_{l+1}(\xi),\ \ P=P_{l+1},\\
P^+=\sum_{j>l}\oplus P_j,\ \ Q=P^+-P,\ \ U=U_{l+1}.
\end{gather*}
As in \cite{Sob}, Section 4,
the estimating of $M(\l; V,\D)$ is reduced
to the problem of eigenvalue estimates for the operators
\begin{gather*}
T_1(\l)=|V|^{1/2}(A-\l I)^{-1}P|V|^{1/2},\\
T_2(\l)=|V|^{1/2}(A-\l I)^{-1}Q|V|^{1/2}.
\end{gather*}
Since $\Vert(A+I)(A-\l I)^{-1}Q\Vert\le C_l$ (actually,
$C_l=O(l)$),
we have
\begin{equation*}
n_+(1,T_2(\l))\le n_+(C^{-1},(|V|^{1/2}(A+I)^{-1}|V|^{1/2})
=
 M(-1; A+C_lV).
\end{equation*}
To the latter quantity the estimate \eqref{w:19}
applies, and we obtain
\begin{equation}\label{w:40}
n_+(1,T_2(\l))\le C'\T(V).
\end{equation}
To the operator $T_1(\l)$ the argument of
\cite{Sob} applies without changes.
Indeed, the
nature of the operator $U$ in our case
is the same as in the case of
periodicity coming from a potential.
This allows to reduce the problem to
estimating the counting function $N(\l; K_V)$
for the operator $K_V$ defined in \eqref{w:15}.
Then using the bound \eqref{w:16},
we arrive at the inequality
$n(1,T_1(\l))\le C(\widetilde J(V)+1)$
which, in combination with \eqref{w:40} leads
to \eqref{w:20}.

(ii) Again, for definiteness,
we prove the result for the non-positive
potentials.
Let $\l = k^2\in I, k >0$.
In the case $l\ge 1$
assume temporarily that $\l\not = (\pi l- \t)^2$
and $\l\not = (\pi l +\t)^2$, so that
$k\in (\pi l -\t, \pi l)\cup(\pi l, \pi l + \t)$.
In the case $l = 0$ assume that $\l \in (0, \t^2)$,
so that $k\in (0, \t)$.
The roots $q_1,\ q_2$ of the
equation \eqref{sp:n3}
are real and distinct, and $q_1 q_2=1$. Let us
label them so that $|q_1|<1<|q_2|$ and
denote $|q_1|=e^{-\s}$, then $\s>0$.
Consider the solutions $y_0$ and $y_1$
(cf. \eqref{sp:sol1} and \eqref{sp:sol})
of the problem
\eqref{sp:eqno} -- \eqref{sp:1}.
Their Wronskian is
$W(y_0, y_1) = y_0 y'_1 - y'_0 y_1
= -b^{1/2}k\sin k\not = 0$,
so that $y_0,\ y_1$
are linearly independent.
On the interval $(n-1,n)$
the function $y_0(t)$ satisfies the inequality
\begin{gather*}
|y_0(t)|\le nb^{1/2}(q_2^{n-2}+q_2^{n-1})
\le 2nb^{1/2}e^{\s t},\ \ n>1;
\notag\\
|y_0(t)|\le |\sin kt|\le kt,\ \ n=1.\notag
\end{gather*}
For $y_1(t)$ we have
\begin{equation*}
|y_1(t)|\le(b^{1/2}+1)q_1^{n-1}
\le q_2(b^{1/2}+1)e^{-\s t},\ \ t>0.
\end{equation*}
Note also that $|q_2|=R|\cos k|+(R^2\cos^2 k-1)^{1/2}\le2R$.
So we see that the inequalities
\begin{equation}\label{w:71}
|y_0(t)|\le cte^{\s t},\ \ |y_1(t)|\le ce^{-\s t},\ \ t>0
\end{equation}
hold uniformly in $\l\in I$.
The solution $y_0$ satisfies the boundary
condition $y_0(0+)=0$.

Given a function $f\in\plainL2(\R_+)$, the
solution of the non-homogeneous equation on $\R_+$:
\begin{equation*}
y''(t)+k^2y(t)=-f(t),\ \ t\not\in\mathbb N;\ \ y(0+)=0
\end{equation*}
satisfying the matching conditions \eqref{sp:1}
for $n\in\mathbb N$, is
given by
\begin{equation*}
y(t)=\int_{\R_+}K(t,s)f(s)ds
\end{equation*}
where
\begin{equation*}
W(y_0,y_1)K(t,s) =
\begin{cases}
y_1(t)y_0(s),\ \ s<t,\\
y_1(s)y_0(t),\ \ t<s.
\end{cases}
\end{equation*}
It follows from \eqref{w:71} that
\begin{multline}
|K(t,s)|\le c^2e^{-\s|t-s|}\min(s,t)
(b^{1/2}k|\sin k|)^{-1}\\
\le C_1 \sqrt{st}(k|\sin k|)^{-1},\ \
C_1 = c^2b^{-1/2}.\label{w:73}
\end{multline}
The operator \eqref{w:6}
(for $\l=k^2$ and $\D=\R_+$) acts as
\begin{equation*}
(T(\l)f)(t)= |V(t)|^{1/2}\int_{\R_+}K(t,s)|V(s)|^{1/2}f(s)ds.
\end{equation*}
Under the assumption
$\widetilde J(V)<\infty$ this operator
belongs to the Hilbert -- Schmidt class.
Indeed, by virtue of \eqref{w:73}
\begin{gather*}
(k\sin k)^2\iint\limits_{\R^2_+}|K(t,s)|^2
|V(t)||V(s)|dtds\\
\le C_1^2
\iint\limits_{\R^2_+}s t |V(t)||V(s)|dtds  =
C_1^2
\biggl(\int\limits_{\R_+} t|V(t)|dt\biggr)^2\le
C_1^2\widetilde J(V)^4.
\end{gather*}
Since $\|T\|\le\|T\|_{HS}$
and $n_+(1,T)=0$ if $\|T\|\le1$, the last
estimate and \eqref{w:7}
imply that $M(\l; A_V) = 0$ if
$\widetilde J(V)^2\le C_1^{-1}k|\sin k|$.
In its turn, this and the estimate
\eqref{w:20} yield the inequality
\eqref{w:201} with
$C' = C\bigl(1+(C_1^{-1}k|\sin k|)^{-1/2}\bigr)$.
By continuity, \eqref{w:201}
extends to the ends of the gap, i.e. to $k=\pi l\pm\t$
($l\ge 1$) or $k = \t$ ($l = 0$).
Since the function
$M(\l; A_V)$ is monotone
in $\l$, the result
for $l=0$ automatically extends to all $\l \le \t$.
\end{proof}

As in \eqref{simplify:eq} one can simplify the estimate
\eqref{w:20} if one assumes that
$|V|$ is decreasing on the interval $\D = (R_1, R_2)$:
\begin{equation}\label{w2:20}
M(\l; A_V, \D)\le
C\biggl[ \biggl(\int_{R_1}^{R_1+1}|V(t)|dt\biggr)^{1/2}
+ \sqrt{6}\int_{\Delta} \sqrt{|V(t)|} dt + 1\biggr].
\end{equation}

\begin{thm}\label{w:asym}
Let the assumptions of \thmref{w:M} be
satisfied and $V\le0$. Then the asymptotics
\begin{equation}\label{w:41}
\lim g^{-1/2}M(\l; gV,\D) = \pi^{-1}
\int_\D\sqrt{|V(t)|}dt,\ \ g\to\infty
\end{equation}
(cf. \eqref{w:17}) holds uniformly in $\l\in\gl_l$.
\end{thm}
\begin{proof}
Like in \cite{Sob}, the problem reduces to the case of a
finite interval
$\D$. In view of \eqref{w:9} we need to study only
the term depending on $g$.
Removal of the matching
conditions inside the interval shifts
the function $N(\l; gV,\D)$ no more than
by $2|\D|+2$ and therefore, does not affect
its asymptotic behaviour.
As a result, we come to the operator of the
Dirichlet problem on a finite
interval for which the asymptotics \eqref{w:41} is well known.
\end{proof}

\subsection{Weyl asymptotics for $M(\l; \BA_{gV})$
and $\widetilde
M(\l; \BA_{gV})$}
The results of this subsection follow immediately from Theorems
3.2 and 3.3.

\begin{thm}\label{w:BA}
Let $\G=\G_b$ and let $V(t)\le0$ be a bounded
measurable function on $\R_+$.
Let $\l$ satisfy \eqref{away:eq}.
Then
\begin{itemize}
\item[(i)]
If
\begin{equation}\label{w:42}
\sum_{k\in\mathbb N}b^k\widetilde J(V_k)<\infty,
\end{equation}
then the Weyl asymptotics holds for the function
$M(\l; \BA_{gV})$
of the operator \eqref{lap:1}:
\begin{equation}\label{w:43}
\lim g^{-1/2}M(\l; \BA_{gV}) = \frac{1}{\pi}
\int_\G\sqrt{|V(|x|)|}dx,\ \ g\to\infty.
\end{equation}
\item[(ii)]
If
\begin{equation}\label{w:44}
\sum_{k\in\mathbb N}\widetilde J(V_k)<\infty,
\end{equation}
then
\begin{equation}\label{w:45}
\lim g^{-1/2}\widetilde M(\l; \BA_{gV}) = \frac{1}{\pi}
\sum_{k=0}^\infty\int_k^\infty
\sqrt{|V(t)|}dt,\ \ g\to\infty.
\end{equation}
The above asymptotic formulae
are uniform in $\l$ on any closed interval $I$
from \eqref{away:eq}.
\end{itemize}
\end{thm}

\begin{proof}
For definiteness, we prove \eqref{w:43}.
The proof of \eqref{w:45}
is the same.

It follows from \eqref{w:2} that
\begin{gather}
g^{-{1/2}}M(\l; \BA_{gV})\label{w:46}\\
= g^{-{1/2}}M(\l; A_{gV})+
(1-b^{-1})\sum_{k\in\mathbb N}
b^kg^{-{1/2}}M(\l; A_{gV_k}).\notag
\end{gather}
By \thmref{w:asym}, for each $k\ge0$
\begin{equation*}
g^{-{1/2}}M(\l; A_{gV_k})
\to\pi^{-1}\int_{\R_+}\sqrt{|V_k(t)|}dt,\ \ g\to\infty.
\end{equation*}
The series \eqref{w:42} dominates
the series
\eqref{w:46} and it follows
from Lebesgue's theorem on the dominated convergence
that
\begin{equation*}
\pi g^{-{1/2}}M(\l; \BA_{gV})
\to \int_{\R_+} \sqrt{|V(t)|}dt + (1-b^{-1})
\sum_{k=1}^\infty b^k\int_{\R_+}
\sqrt{|V_k(t)|}dt,\ \ g\to\infty.
\end{equation*}
This is equivalent to \eqref{w:43}.
Due to \thmref{w:M}
the series \eqref{w:42} converges
uniformly on $I$ from \eqref{away:eq},
and hence the asymptotics
\eqref{w:43} is uniform in $\l\in I$.
\end{proof}

\section{Power-like and exponential potentials}
\label{tree5:sect}

\subsection{Weyl asymptotics}
Here we show how the theorems in the previous section apply to
potentials with a specified rate of decay at infinity.
To have a clear
distinction between the cases of
non-positive and non-negative potentials,
we slightly change our notation:
we denote the potential
by $V=\pm q$ or $V=\sign\! q$ with $\sign = \pm1$,
always assuming that $q\ge  0$.

We are interested in two types of potentials:
power-like and exponential. More precisely, suppose
that $q(t)\le CQ(t)$ where
\begin{equation}\label{Q:eq}
Q(t) = (1+t)^{-2\g}, \ \g>0,\ \
\textup{or} \ \
Q(t) = e^{- 2\vark t},\ \  \vark >0.
\end{equation}
Let us first establish the Weyl type asymptotics for
$M(\l; \BA_{-  gq})$ and $\widetilde M(\l; \BA_{-gq})$.

\begin{thm}\label{powerexp:thm} Let Condition
\eqref{away:eq} be fulfilled.
Suppose that $q(t)\le CQ(t)$.
\begin{itemize}
\item[(i)]
If $Q(t)=(1+t)^{-2\g}$ with $\g>2$, then the asymptotic
formula \eqref{w:45} holds
for $\widetilde M(\l;\BA_{-gq})$;
\item[(ii)]
If $Q(t) = e^{-2\vark t},\ \vark >0$, then
the asymptotic formula \eqref{w:45} for
$\widetilde M(\l;\BA_{-gq})$
holds.
If, in addition, $\vark > \ln b$, then the
asymptotics \eqref{w:43} for $M(\l;\BA_{-gq})$ holds as well.
\end{itemize}
These results are uniform in $\l\in I$
with a closed interval $I$ from \eqref{away:eq}.
\end{thm}

The proof of Theorem \ref{powerexp:thm} is based on two elementary
Lemmas  \ref{power:lem} and \ref{expo:lem} describing individual
counting functions. These Lemmas will be also useful
in the analysis of the non-Weyl behaviour of the
function $\widetilde M$.

Recall that by $q_k,\ k\ge0$ are denoted
the ``shifted'' potentials
$q_k(t)=q(t+k)$, $t>0$.
Remembering the relation \eqref{shift:eq} and a comment
after it, we
sometimes transfer the dependence on
$k$ to the interval $\D$. This is why
some of the estimates below
are stated for intervals $\D$
depending on an additional parameter
$R$, which plays the role of $k$, but is not
supposed to be integer.

\begin{lem}\label{power:lem}
[Power-like potentials]
Suppose that $q\le CQ$ with $Q(t) = (1+t)^{-2\g},\ \g > 1$.
\begin{itemize}
\item[(i)]
Then for any  $R\ge 0$
\begin{equation}\label{known0:eq}
M(\l; \pm  gq, (R, \infty))\le C\bigl(g^{1/2}
(1+R)^{1-\g}+1\bigr),\ \forall k \ge 0,
\end{equation}
uniformly in $\l\in \overline\gl_l$ and $R\ge 0$.
\item[(ii)]
If the condition \eqref{away:eq} is satisfied, then
\begin{equation}\label{known:eq}
M(\l; \pm g q_R)\le C g^{1/2}
(1+R)^{1-\g},\ \forall k \ge 0,
\end{equation}
uniformly in $\l\in I$
with a closed interval $I$ from \eqref{away:eq}.
\end{itemize}
\end{lem}

\begin{lem}\label{expo:lem}[Exponential potentials]
Suppose that $q\le CQ$ with
$Q(t) = e^{-2\vark t},\ \vark > 0$.
\begin{itemize}
\item[(i)] Then for any $R > 0$
\begin{equation}\label{known01:eq}
M(\l; \pm gq, (R, \infty))
\le C\bigl (g^{1/2}e^{-\vark R} +1 \bigr),
\ \forall k\ge 0,
\end{equation}
uniformly in $\l\in \overline{\gl_l}$ and $R\ge 0$.
\item[(ii)]
If
the condition \eqref{away:eq} is satisfied, then
\begin{equation}\label{known1:eq}
M(\l; \pm gq_R)\le C g^{1/2}e^{-\vark R},
\ \forall R\ge 0,
\end{equation}
uniformly in $\l\in I$
with a closed interval $I$ from \eqref{away:eq}.
\end{itemize}
\end{lem}

\begin{proof}[Proofs of Lemmas \ref{power:lem}
and \ref{expo:lem}] Due to the monotonicity of the
function $M(\l; \pm V)$ in $V$,
(see Subsect. \ref{tree4.1:subsect}),
it is sufficient to obtain the estimates for the ``model''
potential $Q$.
For a power-like $Q$, we have
\begin{equation*}
\biggl[\int_R^{R+1}(1+t)^{-2\g}dt\biggr]^{1/2} +
\int_R^\infty (1+t)^{-\g} dt \le  C(1+R)^{1-\g},
\end{equation*}
which implies \eqref{known0:eq} by virtue of \eqref{w2:20}.
Similarly, \eqref{w:201} leads to \eqref{known:eq}.

The proof of Lemma \ref{expo:lem}
is the same.
\end{proof}

\begin{proof}[Proof of Theorem \ref{powerexp:thm}]
According to \eqref{known:eq} (resp. \eqref{known1:eq})
the series \eqref{w:44}
is convergent for $\g >2$ (resp. all $\vark >0$), which ensures
the validity of \eqref{w:45}.

In the exponential case, if $\vark > \ln b$, then the series
\eqref{w:42} is also convergent, which leads to the
asymptotics \eqref{w:43} by Theorem \ref{w:BA}.
\end{proof}

\subsection{Non-Weylian asymptotics}
The rest of the paper is focused on the situations when
the Weyl formula fails, and the asymptotics of
the counting functions \eqref{w:2}--\eqref{w:5}
depends on the behaviour
of the potential at infinity.
We concentrate on bounded potentials
$q$ behaving like $Q$ (see \eqref{Q:eq}) at infinity.
The precise meaning of this phrase will be made clear
later.

As in the previous section, the asymptotics of
$M(\l; \BA_{\pm gq})$
and $\widetilde M(\l; \BA_{\pm gq})$ will be deduced from the
asymptotics of the individual counting functions
$M(\l; \pm gq_k)$, $k\ge 0$ for the operators
for the operators $A_{\pm gq_k}$.
In the case of the power-like potential $q$ the
total number of eigenvalues of $A_{\sign q}$ in
each gap may become infinite. More precisely,
if $q = (1+t)^{-2\g}$, $\sign = -1$ (resp. $\sign = +1$)
and $\g \le 1$ then the eigenvalues accumulate at the
upper (resp. lower)  end of the gap $\gl_l$.
In this connection it
is convenient to introduce the notion
of \textsl{an admissible point} $\l \in \overline\gl_l$.
From now on we fix the number $l\ge 0$ and denote
$\gl_l = (\l_-, \l_+)$.
If $l = 0$, then $\l_- = -\infty$ and $\l_+ = \t^2$.
If $l \ge 1$, then  $\l_\pm = (\pi l \pm \t)^2$.
In the definition below, to avoid unnecessary repetitions,
by $[-\infty, \l_+],\ [-\infty, \l_+)$ we understand the intervals
$(-\infty, \l_+], (-\infty, \l_+)$.

\begin{defn}\label{admis:defn}
Let $Q(t) = (1+t)^{-2\g}$.
Then a point $\l\in [\l_-, \l_+ ]$ is said to
be $\g_0$-admissible, $\g_0 >0$, if
\begin{equation*}
\l\in
\begin{cases}
(\l_-, \l_+],\ \ \g \le \g_0,\ \sign = +1;\\[0.2cm]
[\l_-, \l_+),\ \ \g \le \g_0,\ \sign = -1;\\[0.2cm]
[\l_-, \l_+], \ \ \g > \g_0.
\end{cases}
\end{equation*}
For $Q(t) = e^{-2\vark t}$ any point $\l\in [\l_-, \l_+]$
is said to be $\g_0$-admissible with any $\g_0 >0$.
\end{defn}

Clearly, for any two positive numbers $\g_0, \g_1$,\
$\g_0 < \g_1$, any $\g_1$-admissible $\l$
is automatically $\g_0$-admissible.
For the model potential $Q(t) = (1+t)^{-2\g}$
the number $M(\l; \pm gq_n)$  is finite for all $g >0$
if  $\l$ is $1$-admissible.
For the exponential model potential $Q(t) = e^{-2\vark t}$
the quantity $M(\l; \pm gq_n)$
is finite for all $\l\in[\l_-, \l_+]$.

\subsection{Results for the functions
$\widetilde M(\l; \BA_{\sign gq})$,
$\widetilde N(\l_1,\l_2; \BA_{\sign gq})$}
This subsection contains the results
on the asymptotics of $\widetilde M(\l; \BA_{\pm gq})$ and
$\widetilde N(\l_1, \l_2; \BA_{\pm gq})$.
Their proofs require some technical preparations
which we give in Sections \ref{tree6:sect},
\ref{tree7:sect}.
The proofs are completed in Section \ref{tree8:sect}.
Our results for the functions \eqref{w:2},
\eqref{w:3} require
different techniques and
are much
less complete than those for
$\widetilde M$ and $\widetilde N$;
they are presented in Sect. \ref{tree9:sect}.

Recall that in contrast to the
spectrum of the ``individual'' operator $A_0$
the spectrum of $\BA_0$ contains
eigenvalues $\l_l = (\pi l)^2$ of infinite multiplicity.
Thus, when stating the results
we assume that $\l, \l_1, \l_2$
satisfy \eqref{away:eq} and are
$2$-admissible.
The constants in all the estimates below are
\begin{itemize}
\item
uniform in $\l, \l_1, \l_2$
varying within any closed interval $I$
of $2$-admissible points, satisfying \eqref{away:eq},
\item
independent of the coupling constant $g$.
\end{itemize}
We begin with the power-like potentials.

\begin{thm}\label{main:thm}
Let $q$ satisfy the condition
\begin{equation}\label{asympt:eq}
q(t) = Q(t)\bigl(1+o(1)\bigr),\ t\to \infty,\ \ \
Q(t) = (1+t)^{-2\g},
\end{equation}
and one of the following two conditions
be fulfilled:
\begin{enumerate}
\item
$\g \in (0, 2)$ and $\sign = -1$;
\item
The exponent $\g >0$ is arbitrary and
$\sign = +1$.
\end{enumerate}
Suppose that $\l$ is $2$-admissible and satisfies
\eqref{away:eq}.
Then
\begin{equation}\label{nonweyl:eq}
\lim_{g\to\infty} g^{-\frac{1}{\g}}
\widetilde M(\l; \BA_{\pm g q})
=
\pm\int_0^\infty\int_0^\infty
\bigl[\rho(\l) - \rho\bigl(
\l \mp (s + \s)^{-2\g}
\bigr)\bigr] ds d\s,
\end{equation}
where $\rho$ is the density
of states for the operator $A_0$.
\end{thm}

\begin{rem}A simple change of variables leads to
another expression for the asymptotic coefficient:
\begin{equation*}
\lim g^{-\frac{1}{\g}} \widetilde M(\l; \BA_{\pm gq})
= \pm \int_0^\infty\int_\b^\infty
\bigl[\rho(\l) - \rho\bigl(
\l \mp s^{-2\g}
\bigr)\bigr] ds d\b.
\end{equation*}
In Sect. \ref{tree7:sect} we shall show that the asymptotic
coefficients in the r.h.s. of \eqref{nonweyl:eq} and
that in Theorem \ref{main1:thm} below, are finite.
\end{rem}

Note that in contrast to $\sign = -1$, the
above formula describes the asymptotics of
$\widetilde M(\l; \BA_{\sign gq})$
with $\sign = +1$ for all positive $\g$.
If $\sign = -1 $, then the case $\g = 2$
is critical in the sense that for $\g >2$
the Weyl asymptotics is applicable instead of
\eqref{nonweyl:eq} (cf. Lemma
\ref{power:lem}). We point out however that
for the \textsl{individual} counting function $M(\l; -gq_n)$
the critical case is $\g = 1$ (see Theorem \ref{w:asym}).

To find a formula for $\widetilde M$ in the case $\g = 2$
we need to introduce more restrictions on $q$.

\begin{cond}\label{q:cond}
Let $q\in \plainC1(\R_+)$ be a function such that
\begin{gather*}
c Q(t)\le q(t)\le C Q(t),\ \ \forall t\in \R_+,\\
|q'(t)|\le C Q(t).
\end{gather*}
\end{cond}

Now we are in position to study the critical case:

\begin{thm} \label{g2:thm}
Suppose that $q$ satisfies
\eqref{asympt:eq} with $\g = 2$ and Condition \ref{q:cond}.
Let $\l$ be $2$-admissible and satisfy \eqref{away:eq}.
Then
\begin{equation*}
\lim g^{-1/2} (\ln g)^{-1}\widetilde M( \l; \BA_{-gq})
= (4\pi)^{-1},\ g\to\infty.
\end{equation*}
\end{thm}

The next theorem gives an asymptotic formula for the
number $\widetilde N (\l_1, \l_2)$:

\begin{thm}\label{main1:thm}
Suppose that $q$ satisfies
\eqref{asympt:eq} with some $\g > 0$,
and that
in the case $\sign = -1$, $\a \ge 2$, Condition \ref{q:cond}
is also fulfilled.
Let $\l_1, \l_2$ be $2$-admissible and satisfy
\eqref{away:eq}. Then
\begin{multline}\label{nonweyl1:eq}
\lim g^{-\frac{1}{\g}} \widetilde N (\l_1, \l_2, \BA_{\pm gq})\\
= \int_0^\infty
\int_0^\infty \bigl[
\rho\bigl( \l_2 \mp (t + \s)^{-2\g} \bigr)
- \rho\bigl( \l_1 \mp (t + \s)^{-2\g} \bigr)\bigr]dt d\s,
\end{multline}
as $g\to\infty$.
\end{thm}

We point out that the asymptotics of
$\widetilde N(\l_1, \l_2)$ is described
by the density of states $\rho(\l)$ for \textsl{all} $\g >0$.
Under the conditions of Theorem \ref{main:thm}
the asymptotics \eqref{nonweyl1:eq} can be immediately deduced
from \eqref{nonweyl:eq} with the help of \eqref{t:eq}.
On the contrary, for $\a > 2$ and $\sign = -1$
the behaviour of $\widetilde N$ can not be inferred from the
asymptotics of $\widetilde M(\l)$ which is given by the Weyl
term, see Lemma \ref{power:lem}.

Let us proceed to the exponential potentials.
From Lemma \ref{expo:lem}
we know  that for the case $q\le CQ$,
$Q(t) = e^{-2\vark t}$, $\sign = -1$, the
asymptotics of $\widetilde M(\l; \BA_{\sign q})$
is described by the Weyl formula \eqref{w:45}.
The next theorem
gives an answer in the case $\sign = +1$.
Below $g_0 > e$ is a constant.

\begin{thm}\label{exp+:thm}
Suppose that
\begin{equation}\label{twoside:eq}
cQ(t)\le q(t)\le CQ(t),\ \forall t\ge R_0,
\end{equation}
with $Q(t) = e^{-2\vark  t}, \vark >0$ and some $R_0\ge 0$.
Let $\l, \l_1, \l_2$ be arbitrary numbers satisfying
\eqref{away:eq}.
Then
\begin{equation}\label{tildemexp:eq}
\widetilde M(\l; \BA_{gq}) = \frac{1}{8\vark^2}
\rho(\l)(\ln g)^2 + O(\ln g),\  g \ge g_0,
\end{equation}
and
\begin{equation}\label{tildenexp+:eq}
\widetilde N(\l_1, \l_2; \BA_{gq})
\le C\ln g,\ g\ge g_0.
\end{equation}
\end{thm}

The next result complements the Weyl formula \eqref{w:45}
by providing an estimate for the function
$\widetilde N(\l_1, \l_2; \BA_{-gq})$:

\begin{thm}\label{exp-:thm}
Suppose that $q$ fulfills Condition \ref{q:cond}
with $Q(t) = e^{-2\vark t}$.
Let $\l_1, \l_2$ be arbitrary numbers satisfying \eqref{away:eq}.
Then
\begin{equation}\label{tildenexp-:eq}
\widetilde N(\l_1, \l_2; \BA_{-gq})\le C(\ln g)^2,\
g\ge g_0.
\end{equation}
\end{thm}

\section{Individual estimates and Weyl asymptotics
with a remainder}\label{tree6:sect}

\subsection{Individual estimates}

Here we obtain further estimates for individual
counting functions $M(\l; \pm gq, \D)$.
Although our ultimate objective is to establish asymptotic
formulae for the counting functions
$M(\l; \pm gq_k)$ with integer
non-negative $k$'s, most of the results in this section
are uniform with respect to a wide class of
potentials, including the shifted potentials $q_R, R \ge 0$.

Unless stated otherwise, in this section we always assume
that the points $\l, \l_1, \l_2\in [\l_-, \l_+]$
are $1$-admissible. The constants in all the
estimates obtained below are
\begin{itemize}
\item
uniform in $\l, \l_1, \l_2$
varying within any closed interval $I\subset [\l_-, \l_+]$
of $1$-admissible points;
\item
independent of the coupling constant $g$.
\end{itemize}
Whenever possible we treat the power-like
and exponential potentials
simultaneously. It is convenient to use the notation
\begin{equation}\label{alpha:eq}
\a =
\begin{cases}
g^{\frac{1}{2\g}},\ \
\textup{if}\ \  Q(t) = (1+t)^{-2\g};\\[0.2cm]
\dfrac{1}{2\vark}\ln g,\ \
\textup{if}\ \  Q(t) = e^{-2\vark t}.
\end{cases}
\end{equation}
We always assume that $g\ge g_0> e$,
so $\a\ge \a_0$ with
some $\a_0> 0$ in both cases.

The following simple Lemma will be
repeatedly used:

\begin{lem}\label{decoupl:lem}
Suppose that $q(t)\le CQ(t)$.
\begin{itemize}
\item[(i)]
If $Q(t) = (1+t)^{-2\g}$,
then for all $R\ge 1$
\begin{equation}\label{decoupl1:eq}
\bigl|M(\l; \pm gq, \R_+)
- M(\l; \pm gq, (0, R\a))\bigr| \le C'\a.
\end{equation}
with a constant $C'$ depending only on $C$.
Moreover,
\begin{equation}\label{decoupl2:eq}
\limsup_{R\to\infty}
\sup_{g\ge g_0}
\a^{-1}\bigl|M(\l; \pm gq, \R_+)
- M(\l; \pm gq, (0, R\a))\bigr| = 0.
\end{equation}
\item[(ii)]
If $Q(t) = e^{-2\vark t}$,
then
\begin{equation}\label{decoupl3:eq}
\bigl|M(\l; \pm gq, \R_+)
- M(\l; \pm gq, (0, R\a))\bigr| \le C',
\end{equation}
with a constant $C'$ depending only on $C$.
\end{itemize}

\end{lem}

\begin{proof}
In view of the decoupling principle \eqref{w:13}
\begin{equation*}
\bigl| M(\l; \pm gq, \R_+)
- M(\l; \pm gq, (0, R\a))\bigr| \le 2
+ M(\l; \pm qq, (R\a, \infty) ).
\end{equation*}
If $Q(t) = (1+t)^{-2\g}$ and $\g > 1$, then by
\eqref{known0:eq}
the last term in the r.h.s.
does not exceed
\begin{equation*}
g^{1/2}( \a R + 1)^{1-\g} + 1 \le \a R^{1-\g} + 1.
\end{equation*}
This implies \eqref{decoupl1:eq}
and \eqref{decoupl2:eq}.
If $\g\le 1$, then  $\l < \l_+$ (for $\sign= -1$) or
$\l > \l_-$ ( for $\sign = +1$). Thus
one can choose $\hat R$ so as to ensure that
$M(\l; \pm qq, (\hat R\a, \infty) ) = 0$, since
$|gq(t)|\le g \hat R^{-2\g} \a^{-2\g}
= \hat R^{-2\g}$ for
all $t\ge \hat R\a$. Now \eqref{decoupl2:eq} follows.
To show \eqref{decoupl1:eq} use the
decoupling principle and
\eqref{w2:20} to conclude that
\begin{align*}
M(\l; \pm gq, (R\a, \infty) )
\le &\ 2+ M(\l; \pm gq, (R\a, \hat R\a) )\\[0.2cm]
\le &\ C \biggl[\a \biggl[\int_R^{R+1}
|t|^{-2\g} dt\bigg]^{1/2}
+ \a\int_R^{\hat R} |t|^{-\g} dt + 1\biggr] + 2\\[0.2cm]
\le &\ C(1+\a).
\end{align*}
For the case $Q(t) = e^{-2\vark t}$,
using \eqref{known01:eq}
we obtain the estimate:
\begin{equation*}
M(\l; \pm gq, (R\a, \infty)) \le
C\bigl(g^{1/2} e^{-\vark\a R}  + 1\bigr)
\le C( g^{1/2-R/2} + 1)
\le C
\end{equation*}
for all $R\ge 1$, as required.
\end{proof}

\begin{lem}\label{est:lem}
Let $q(t)\le CQ(t)$.
\begin{itemize}
\item[(i)]
If  $Q(t) = (1+t)^{-2\g}$,
then
\begin{equation}\label{estpm:eq}
M(\l; - gq)\le
\begin{cases}
Cg^{1/2},\ \g > 1,\\[0.2cm]
Cg^{1/2} \ln g,\ \g = 1,\\[0.2cm]
C\a,\ \g < 1.
\end{cases}
\end{equation}
\item[(ii)]
Let either $Q(t) = (1+t)^{-2\g}$
with arbitrary $\g >0$,
or $Q(t) = e^{-2\vark t}$.
Then
\begin{equation}\label{est+:eq}
M(\l; gq)\le C \a,
\end{equation}
and
\begin{equation}\label{est:eq}
N(\l_1, \l_2; gq)\le C\a.
\end{equation}
\end{itemize}
\end{lem}

\begin{proof}
(i) The estimate \eqref{estpm:eq} for $\g >1$
follows from \eqref{known0:eq}.
If $\g \le 1$, then by  \eqref{w2:20},
\begin{align*}
M(\l; - gq, (0, \a))
\le &\ Cg^{1/2}
\biggl[\int_0^1 (1+t)^{-2\g} dt\biggr]^{1/2}
+ Cg^{1/2} \int_0^{\a} t^{-\g} dt + C\\
\le &\
\begin{cases}
Cg^{1/2} \ln \a,\ \g = 1,\\[0.2cm]
C\a,\ \g < 1.
\end{cases}
\end{align*}
By virtue of Lemma \ref{decoupl:lem} this
leads to \eqref{estpm:eq}.

(ii) It is clear from \eqref{w:9}
that
\begin{equation*}
M(\l; gq, (0, \a))\le N(\l; 0, (0, \a) )\le C\a.
\end{equation*}
Now Lemma \ref{decoupl:lem} gives \eqref{est+:eq}.
The estimate \eqref{est:eq} follows from \eqref{est+:eq}
and \eqref{w:11}.
\end{proof}

From these Lemmas we can immediately deduce the
asymptotics of $M(\l; g q)$
with an exponential $q$.

\begin{thm}\label{plus:thm}
Let $q$ be a bounded function  satisfying
\eqref{twoside:eq} with $Q(t) = e^{-2\vark t}$.
Then
\begin{equation}\label{mexp:eq}
\bigl| M(\l; gq) - \rho(\l)\a\bigr|\le C,\ \forall \a\ge 1,
\end{equation}
and
\begin{equation}\label{nexp:eq}
N(\l_1, \l_2; gq)\le C,\ \forall \a \ge 1.
\end{equation}
\end{thm}

\begin{proof}
The bound for $N(\l_1, \l_2; gq)$ immediately
follows from \eqref{mexp:eq} by \eqref{t:eq},
since $\rho(\l_1) = \rho(\l_2)$ for $\l_1, \l_2\in[\l_-, \l_+]$.

By the decoupling principle \eqref{decoupl3:eq}
it suffices to study the
counting functions $M(\l; gq, \Delta)$
with $\Delta = (0, \a)$.
The result for the unperturbed function $N(\l; 0, \Delta)$
immediately follows from Theorem \ref{density:thm}:
\begin{equation}\label{free:eq}
N(\l; 0, \Delta) = \rho(\l)\a + O(1),\ \a\to\infty.
\end{equation}
To handle
the perturbed function
split the interval $\Delta$ as follows:
\begin{gather*}
\Delta = \D_0\cup \Delta_1 \cup \Delta_2,\\
\D_0 = (0, R_0],\
\Delta_1 = (R_0, R_1\a],\
\Delta_2 = (R_1 \a, \a),
\end{gather*}
where $R_0>0$ is defined in \eqref{twoside:eq}.
The number $R_1>0$ is found from the requirement
\begin{equation*}
gq(t) \ge \l - \l_0,\  \forall t\in (R_0, R_1\a],
\end{equation*}
where $\l_0 = \t^2 = \inf\s(A_0)$.
This implies that for $R_1$ we can take
\begin{equation}\label{r1:eq}
R_1 = 1 - \frac{C'}{2\vark\a},
\end{equation}
with a sufficiently large $C' = C'(\l) >0$.
It follows from \eqref{r1:eq}
and Theorem \ref{density:thm} that
\begin{equation*}
N(\l; gq, \Delta_2)
\le N(\l; 0, \Delta_2)\le (1 - R_1)\a + C'
\le C'',\ \forall \a\ge 1,
\end{equation*}
Since $N(\l; gq, \D_0)\le N(\l; 0, \D_0)\le C$ and
$N(\l; gq, \Delta_1) = 0$, by the decoupling
principle \eqref{w:12} we have
\begin{equation*}
N(\l; gq, \Delta)\le C,\ \forall \a\ge 1.
\end{equation*}
Now it follows from \eqref{free:eq}
and \eqref{w:9} that
\begin{equation*}
\bigl|M(\l; gq, \Delta) - \rho(\l)\a\bigr| \le C,\ \a\ge 1,
\end{equation*}
which implies \eqref{mexp:eq}.
\end{proof}

\subsection{Individual asymptotics of
the Weyl type with a remainder}
\label{weyl:subsect}
Even if a potential $V\le 0$ satisfies the conditions
of Theorem \ref{w:asym}, the formula \eqref{w:41}
fails to provide an asymptotics for
$N(\l_1, \l_2; V)$ as the leading term in \eqref{w:41}
does not depend on $\l$.
Below we establish, under certain conditions on $q$, a
Weyl-type asymptotics for $M(\l; -gq, \D)$ with a remainder, which
allows us to obtain bounds on the growth of $N(\l_1, \l_2; -gq)$
as $g\to\infty$.

We begin with an asymptotic formula for the Schr\"odinger
operator $ - d^2/dt^2 - q$
on a bounded interval $\D\subset\R_+$
with the Dirichlet boundary conditions
but \textsl{without} any matching
conditions. Denote the counting function of this operator
by $\#(\l; -q, \D)$.
The next theorem is a minor modification
of a similar statement from \cite{Sob}:

\begin{thm}\label{remainder:thm}
Let $q\in \plainC1([0, \infty))$ be a non-negative function,
and let $\D = (R_1, R_2)$ with
$0\le R_1 \le R_2 < \infty$.
Then for any  $\l\in\R$
one has
\begin{multline}
\biggl|\#(\l; -q; \Delta)  - \frac{1}{\pi}
\int_{\D} \sqrt{q(t)} \ dt\biggr|\\[0.2cm]
\le \int_\D \frac{|q'(t)|}{4\pi (q(t)+|\l|)}\ dt +
\frac{3\sqrt{|\l|+1}}{\pi}\ |\D|  + 1,\label{weyl10:eq}
\end{multline}
where the constant $C$ does not depend on $\D, g$ and
is uniform in $\l$ on a compact interval.
\end{thm}

\begin{proof}
Assume without loss of generality that $R_1 = 0$.
The idea is to use the fact that the
number $\#(\l) = \#(\l; -q, \D)$ equals
the number of roots of the solution $u$ of the equation
\begin{equation}\label{equation:eq}
-u'' - q u = \l u,\ u(0) = 0,\  u'(0) = 1,
\end{equation}
lying strictly inside $\D$. To find the number
of roots, represent $u$ in the polar form:
\begin{equation}\label{polar:eq}
u(t) = \b(t) \sin\xi(t),\ \
u'(t) = \b(t)f(t)\cos\xi(t),\ f = \sqrt{q + \l_0},
\end{equation}
with a $\l_0>0$.
The equation \eqref{equation:eq} and the above equalities
define the real-valued amplitude $\b$
and the phase $\xi$ uniquely
under the assumption that $\xi$ is continuous.
Substituting \eqref{polar:eq}
in the equation \eqref{equation:eq}, one obtains the following
non-linear equation for $\xi$:
\begin{equation}\label{xi:eq}
\xi' = f + \frac{f'}{2f}\sin(2\xi)
+ \frac{\l-\l_0}{f} \sin^2\xi, \ \ \xi(0) = 0,
\end{equation}
and a linear equation for $\b$:
\begin{equation*}
\b' = - \frac{1}{f}\bigl((\l-\l_0)\sin\xi
+ f' \cos\xi\bigr)\cos\xi \ \b ,
 \ \b(0) = \frac{1}{f(0)}.
\end{equation*}
Since $\b$ never vanishes, the number of roots of $u$
equals the number of points $t\in\D$
where $\xi(t) = 0(\modulo \pi)$.
From \eqref{xi:eq} it is clear that $\xi'(t) = f(t) > 0$
for those $t$, so that
\begin{equation*}
\pi^{-1}\xi(R_2) - 1 \le \#(\l)\le \pi^{-1}\xi(R_2).
\end{equation*}
Therefore \eqref{xi:eq} implies that
\begin{equation*}
\biggl|\#(\l) - \frac{1}{\pi}\int_{\D} f dt\biggr|
\le \int_{\D}\frac{|f'|}{2\pi f} dt
+ \int_{\D}\frac{|\l-\l_0|}{\pi f} dt + 1.
\end{equation*}
Since $f\ge \sqrt{\l_0}$,
$f' = q' (2f)^{-1}$, and
$\sqrt{q+\l_0} - \sqrt{q}\le \sqrt{\l_0}$,
this leads to
\begin{equation*}
\biggl|\#(\l) - \frac{1}{\pi}\int_{\D} \sqrt{q}\  dt\biggr|
\le \int_{\D}\frac{|q'|}{4\pi (q + \l_0)} dt
+ \int_{\D}\frac{|\l| + 2 \l_0}{\pi \sqrt{\l_0}} dt + 1.
\end{equation*}
It remains to take $\l_0 = |\l|+1$.
\end{proof}

\begin{thm}\label{weyl:thm}
Suppose
that $q$ satisfies Condition \ref{q:cond} and
let $\D = (R_1, R_2)$ with
$0\le R_1 \le R_2 < \infty$.
Then for any $\l\in\R$ one has
\begin{equation}\label{weyl11:eq}
\biggl|N(\l; -gq; \Delta)  - \frac{\sqrt g}{\pi}
\int_{\D} \sqrt{q(t)}\   dt\biggr|
\le C(|\D| + 1),
\end{equation}
for all $g\ge g_0$,
where the constant $C$ does not depend on $\D, g$ and
is uniform in $\l$ on a compact interval.
\end{thm}

\begin{proof}
Let
us split $\Delta$ as follows:
\begin{gather*}
\Delta = \cup_{k=l}^m \Delta_k,\ l = [R_1], m = [R_2],\\
\Delta_l = (R_1, l+1],\ \Delta_m = (m, R_2 ),\\
\Delta_k = (k, k+1],\ k = l+1, \dots, m-1.
\end{gather*}
Then by the decoupling principle
\begin{equation*}
\bigl| N(\l; \Delta) - \sum_k N(\l; \Delta_k)\bigr|
\le 2(m-l)\le 2(|\D| + 2).
\end{equation*}
To study each $\D_k$ we use Theorem \ref{remainder:thm}.
Namely, since $|q'|\le Q$ and $q\ge cQ$,
the estimate \eqref{weyl10:eq} yields:
\begin{equation*}
\biggl|N(\l; \Delta_k)
- \frac{\sqrt{g}}{\pi}\int_{\D_k} \sqrt{q(t)}\  dt\biggr|
\le C
\end{equation*}
with a constant independent of $k$.
Adding up these inequalities over
$k = l, l+1, \dots, m$, we arrive at
\eqref{weyl11:eq}.
\end{proof}

Let us derive from this theorem  the asymptotics for
the counting function $M(\l; -gq, \D)$ with
$\D = (k, \infty)$, $k >0$.

\begin{thm}\label{weyl1:thm}
Suppose that Condition \ref{q:cond} is satisfied
and let $\D = (k, \infty)$, $k\in \mathbb N \cup \{0\}$.
Then
\begin{equation}\label{weyl1:eq}
\biggl|M(\l; -gq, \D)  - \frac{\sqrt g}{\pi}
\int_k^{k+\a} \sqrt{q(t)} \ dt\biggr| \le C\a,
\ \forall g\ge g_0,
\end{equation}
uniformly in $k$,
and
\begin{equation}\label{n:eq}
N(\l_1, \l_2; -gq, \D)\le C\a, \ \forall g\ge g_0.
\end{equation}
If $Q(t) = (1+t)^{-2\g}, \g >1$ or
$Q(t) = e^{-2\vark t}$, then
\begin{equation}\label{weyl2:eq}
\biggl|M(\l; -gq, \D)  - \frac{\sqrt g}{\pi}
\int_k^\infty \sqrt{q(t)} \ dt\biggr| \le C\a,
\ \forall g\ge g_0.
\end{equation}
\end{thm}

\begin{proof} The estimate \eqref{n:eq}
follows from \eqref{weyl1:eq}
by virtue of \eqref{w:11}.

Let us prove \eqref{weyl1:eq}.
Let $\Delta_1 = (k, k+ \a]$ and
$\Delta_2 = (k+\a, \infty)$.
In view of Lemma \ref{decoupl:lem} and the relation
\eqref{shift:eq}
it suffices to show that the distribution function
$M(\l; -gq, \Delta_1)$ satisfies \eqref{weyl1:eq}.
This is a direct consequence of Theorem \ref{weyl:thm}
and the identity \eqref{w:9}.

The formula \eqref{weyl2:eq} follows from
\eqref{weyl1:eq} in view of the inequality
\begin{equation*}
\sqrt g \int_{k+\a}^\infty \sqrt {q(t)} dt \le C \a.
\end{equation*}
\end{proof}

\section{Individual asymptotics. Power-like potentials}
\label{tree7:sect}

In this section we study the individual counting function
$M(\l; \sign gq_k)$ for
$q$ satisfying \eqref{asympt:eq} with
$Q(t) = (1+t)^{-2\g}$, where $\g >0$ is arbitrary for both cases
$\sign = \pm 1$.
In contrast to the previous section
here we focus on the asymptotics  of this function
under the assumption that $g$ and $k$ tend to infinity in a
coordinated way (see Lemma \ref{individ:lem} below).
We shall use the following notation:
\begin{equation}\label{ab:eq}
\a = g^{\frac{1}{2\g}},\ \b = \b_k(g) = \frac{k+1}{\a}.
\end{equation}
Emphasise again that the main difference with the
asymptotics obtained
in Theorem \ref{w:asym} is that now
it is determined by the density of
states for the unperturbed
operator $A_0$.

We begin with the study of asymptotic coefficients.

\subsection{Asymptotic coefficients}
Introduce the asymptotic coefficients for
$M(\l; \pm gq_k)$:
\begin{equation}\label{fpm:eq}
F_\pm(\s, \l)
=  \pm\int_0^\infty \bigl[\rho(\l) -
\rho\bigl( \l \mp (s + \s)^{-2\g}
\bigr)\bigr]ds
\end{equation}
and for $N(\l_1, \l_2; \pm gq_k)$:
\begin{multline}\label{gpm:eq}
G_\pm(\s, \l_1, \l_2) =
\pm \bigl(F_\pm(\s, \l_1)
- F_\pm(\s, \l_2)\bigr)\\[0.2cm]
= \int_0^\infty \bigl[
\rho\bigl( \l_2 \mp (s + \s)^{-2\g} \bigr)
- \rho\bigl( \l_1 \mp (s + \s)^{-2\g} \bigr)\bigr]ds.
\end{multline}
It is clear that $F_\pm \ge 0$ and  $G_\pm \ge 0$
if $\l_1 \le \l_2$.
Some other useful properties of $F_\pm, G_\pm$
are collected in the next Lemma:

\begin{lem}\label{f:lem}
Let $\l, \l_1, \l_2\in [\l_-, \l_+]$ be $1$-admissible
numbers. Then the integral $F_\pm(\s, \l)$ is finite for all
$\s >0$. Moreover,
\begin{itemize}
\item[(i)]
If $\g > 1$, then
\begin{gather*}
F_-(\s, \l)\le C \s^{1-\g},\\[0.2cm]
F_+(\s, \l)\le C(1+\s)^{1-\g},
\end{gather*}
for all $\s>0$.
\item[(ii)]
If $\g\le 1$, then $F_+(\s, \l)\le C$ and
\begin{equation*}
F_-(\s, \l)\le
\begin{cases}
C,\ \g <1,\\
C\ln (\s^{-1} +1),\  \g=1.
\end{cases}
\end{equation*}
\item[(iii)]
If $\sign = +1$ and $\l > \l_-$ or
$\sign = -1$ and $\l < \l_+$, then
$F_{\pm}(\s, \l) = 0$ for all
$\s\ge \s_{\pm}(\l) = \dc_{\pm}^{-\frac{1}{2\g}}$ with
$\dc_{\pm} = |\l - \l_{\mp}|$, and
\begin{equation*}
F_{\pm} (\s, \l) \ge C(\s_{\pm} - \s)^{\frac{3}{2}},\ \s\le \s_\pm.
\end{equation*}
\item[(iv)]
The integral
\begin{equation}\label{f1:eq}
F_\pm^R(\s ,\l) = \pm \int_0^R
\bigl[ \rho(\l) - \rho\bigl(
\l \mp (s + \s)^{-2\g}\bigr)\bigr]ds,\  R >0,
\end{equation}
tends to $F_\pm(\s, \l)$ as $R\to\infty$ uniformly in $\s >0$.
\item[(v)]
For all $\g>0$ one has
$G_\pm(\s, \l_1, \l_2)\le C$,\ $\forall \s >0$.
\end{itemize}

\end{lem}

\begin{proof}
(i)
By \eqref{holder1:eq},
the integrand in the definition \eqref{fpm:eq} does not
exceed $C(s+\s)^{-\g}$ for $\sign = -1$ and
$\min\{\rho(\l), C(s+\s)^{-\g}\}$ for $\sign = +1$.
The required estimates follow
immediately.

(ii) Let  $\g \le 1$.
Let first $\sign = +1$, so that $\l\in (\l_-, \l_+]$.
Define
$R >0$ to be the number such that $\l-R^{-2\g} = \l_-$.
Consequently, $\l-(s+\s)^{-2\g} \ge \l_-$ for
$s+\s\ge R$, and
hence
\begin{equation*}
\rho(\l) - \rho\bigl(\l-(s+\s)^{-2\g}\bigr)= 0,\
\forall s: s+\s\ge R.
\end{equation*}
This implies that
\begin{equation*}
F_+(\s, \l)\le \int_0^R \rho(\l) ds
= R\rho(\l)\le C'(\l).
\end{equation*}
Consider now the case $\sign = -1$,
so that $\l\in [\l_-, \l_+)$.
Let $R >0$ be the number such that $\l+R^{-2\g} = \l_+$.
Consequently, $\l+(s+\s)^{-2\g} \le \l_+$
for $s+\s\ge R$, and
hence
\begin{equation*}
\rho\bigl(\l+(s+\s)^{-2\g}\bigr) - \rho(\l) = 0,\
\forall s: s+\s\ge R.
\end{equation*}
If $\s\ge R$, then $F_-(\s, \l) = 0$. If $\s < R$, then
\begin{equation*}
F_-(\s, \l)\le C\int_0^{R - \s} (\s+s)^{-\g} ds
\le
\begin{cases}
C(R),\ \g < 1,\\[0.2cm]
C(R) \ln(\s^{-1} + 1),\ \g=1.
\end{cases}
\end{equation*}

(iii)
For brevity consider only the case
$\sign = -1$, so that $\l < \l_+$.
For $s \ge \s_- - \s, \ \s_- = \s_-(\l)$, the integrand in
\eqref{fpm:eq} equals zero.
Besides, as $\rho(\l) = \rho(\l_+)$,
in view of \eqref{effmass1:eq} we have
\begin{gather*}
\rho\bigl(\l+(s+\s)^{-2\g}\bigr) - \rho(\l)
\ge c((s+\s)^{-2\g} - \dc_-)^{\frac{1}{2}}
\ge c'(\s_- -\s - s)^{\frac{1}{2}},\\
\forall s\in (0, \s_- - \s).
\end{gather*}
Integrating this inequality in $s$, we obtain the
required lower bound for $F_-(\s, \l)$.
The analogous bound for $F_+(\s, \l)$ is
obtained in the same way.

(iv) It suffices to notice that for any
$1$-admissible $\l$ and any $R >0$ one has
\begin{equation*}
\pm \int_R^\infty \biggl[
\rho(\l)
- \rho\bigl(\l \mp (\s+s)^{-2\g}\bigr)\biggr] ds
\le \pm \int_R^\infty \biggl[
\rho(\l)
- \rho\bigl(\l \mp s^{-2\g}\bigr)\biggr] ds
\end{equation*}
for all $\s \ge 0$.

(v) Arguing as on the previous step,
it suffices to prove that the integral of the form
\eqref{gpm:eq} over a finite interval
is bounded uniformly in $\s >0$.
By \eqref{holder1:eq}
\begin{equation*}
\bigl|
\rho\bigl( \l_2 \mp (s + \s)^{-2\g} \bigr)
- \rho\bigl( \l_1 \mp (s + \s)^{-2\g} \bigr)\bigr|
\le C|\l_1-\l_2|^{1/2},
\end{equation*}
which provides the required boundedness.
\end{proof}

When studying the sum of the
counting functions, we shall
need some properties of the sum of asymptotic
coefficients $F_\pm(\b_k, \l)$:

\begin{lem}\label{sumf:lem}
\begin{itemize}
\item[(i)]
Suppose that $\l$ is a $1$-admissible number and $\d=\d(\a)$
is a bounded function such that
\begin{equation}\label{small:eq}
\d\le 1,\  \a\d^{2\g+1} \to\infty,\ \  {\textup{as}}\ \ \a\to\infty.
\end{equation}
Then for any fixed $A\ge \sup_\a\d(\a)$ one has
\begin{equation}\label{sumf:eq}
\lim \biggl[\a^{-1}\sum_{[\d\a]\le k\le [A\a]}
F_\pm(\b_k, \l) - \int_\d^A F_\pm(\s,\l) d\s\biggr] = 0,\ g\to\infty.
\end{equation}
\item[(ii)]
If $\l, \l_1, \l_2$ are $2$-admissible, then
the integrals
\begin{equation*}
\int_0^\infty F_+(\s,\l)d\s,\
\int_0^\infty G_\pm(\s,\l_1, \l_2)d\s
\end{equation*}
are finite.
If in addition $\g \in (0, 2)$, then the integral
\begin{equation*}
\int_0^\infty F_-(\s,\l)d\s
\end{equation*}
is also finite.
\end{itemize}
\end{lem}

\begin{proof} For brevity we omit $\l$ from the notation of
$F_\pm$.

(i) Let $\s\in (k, k+1]$ be an arbitrary number,
and let
\begin{equation*}
\l_1 = \l\pm(s+\b_k)^{-2\g}, \
\l_2 = \l\pm(s+\s/\a)^{-2\g},\ t >0.
\end{equation*}
Observe that
\begin{equation*}
|\l_1-\l_2| \le  2\g \d^{-2\g-1}\a^{-1}.
\end{equation*}
Now it follows from \eqref{holder1:eq} that
\begin{equation*}
\bigl|\rho(\l_1) - \rho(\l_2)
\bigr|
\le CE_1(\a, \d),\ \
E_1(\a, \d) = \a^{-1/2} \d^{-\g -1/2}.
\end{equation*}
Thus, by definition \eqref{f1:eq}, for each $R >0$ one has
\begin{equation*}
\biggl|F_\pm^R(\b_k)
- \int_{k}^{k+1} F_\pm^R(\s/\a) d\s\biggr|
\le C R E_1(\a, \d),
\end{equation*}
or, changing the variable under the integral,
\begin{equation*}
\biggl|F_\pm^R(\b_k)
- \a\int_{k/\a}^{(k+1)/\a} F_\pm^R(\s) d\s\biggr|
\le C R E_1(\a, \d).
\end{equation*}
Therefore
\begin{multline*}
\left|\a^{-1}\sum_{k=[\d \a]}^{[A\a]} F_\pm^R(\b_k)
- \int_\d^A F_\pm^R(\s) d\s \right|\\
\le CAR E_1(\a, \d)
+ \int_{[\d\a]\a^{-1}}^\d F_\pm^R(\s)d\s
+ \int_A^{([A\a]+1)\a^{-1}}
F_\pm^R(\s) d\s.
\end{multline*}
Clearly, the first term tends to zero under the conditions
\eqref{small:eq}.
The last two integrals tend to zero as $\a\to\infty$ by
Lemma \ref{f:lem}(i), (ii).
Consequently, for each $R>0$
\begin{align*}
\limsup \biggl|\a^{-1}\sum_{k=[\d\a]}^{[A\a]}&\
F_\pm(\b_k)  -  \int_{\d}^A F_\pm(\s) d\s
\biggr|\\
\le &\ A\limsup \max_{[\d\a] \le k\le [A\a]}
\bigl| F_\pm^R(\b_k) - F_\pm(\b_k)\bigr|\\
&\ + \int_{\d}^A \bigl|F_\pm^R(\s) - F_\pm(\s)\bigr| d\s,
\end{align*}
where $\limsup$ is taken under the conditions \eqref{small:eq}.
Recall that by Lemma \ref{f:lem}
$F_\pm^R(\s)$ converges to $F_\pm(\s)$ as $R\to\infty$
uniformly in $\s>0$. Thus the r.h.s. of the above
inequality vanishes as $R\to\infty$.
This proves \eqref{sumf:eq}.

(ii) By Lemma \ref{f:lem}(i), (v), and also by definition
\eqref{gpm:eq}, the functions
$F_+$ and $G_\pm$ are integrable in $\s$ for $\g >2$.

If $\g < 2$, then  $F_\pm, G_\pm$ have compact support
due to Lemma \ref{f:lem}(iii).
They are also integrable near $\s = 0$ by virtue
of Lemma \ref{f:lem}(i).

If $\g = 2$, then the same applies to $F_+$ and $G_\pm$
again by Lemma \ref{f:lem}(i), (iii), (v).
\end{proof}

Lemma \ref{sumf:lem} guarantees that
the asymptotic coefficients in Theorems
\ref{main:thm} and \ref{main1:thm} are finite.

\subsection{Asymptotics of $M(\l; \pm gq_k)$}

\begin{lem}\label{individ:lem}
Let $q$ satisfy \eqref{asympt:eq},
and let $\l, \l_1, \l_2\in[\l_-, \l_+]$
be $1$-admissible. Then for any  function $\d = \d(\a)$
satisfying \eqref{small:eq} and any fixed $A > \sup_\a\d(\a)$,
one has
\begin{gather}
\lim \max_{\d\le \b_k \le A}
\bigl|\a^{-1} M(\l; \pm g q_k)
- F_\pm\bigl(\b_k(g), \l\bigr)\bigr| = 0,
\label{final:eq}\\[0.2cm]
\lim \max_{\d\le \b_k \le A}
\bigl|\a^{-1} N(\l_1, \l_2; \pm g q_k)
- G_\pm\bigl(\b_k(g), \l_1, \l_2\bigr)\bigr| = 0,
\label{final1:eq}
\end{gather}
as $g\to \infty$.
\end{lem}

\begin{proof}
Note without further ado,
that \eqref{final1:eq} is a direct consequence of
\eqref{final:eq} in view of \eqref{t:eq}.

Let us concentrate on the proof of \eqref{final:eq}.
By \eqref{decoupl2:eq}
it suffices to establish the required
asymptotic formula for the counting function
$M(\l; \pm gq_k, \Delta)$ with
$\Delta(g) = (0, R\a]$ and afterwards take $R$ to infinity.
The proof of this fact is an adaptation of the
corresponding argument from \cite{Sob}.

Suppose first that $q(t) = Q(t)$.
Let us split $(0, R]$ into $L$ identical subintervals
$(s_{j-1}, s_j],\ j = 1, 2, \dots, L$, so that
$s_0 = 0, s_L = R$ and $s_{j+1} - s_j = RL^{-1}$,
and denote
\begin{equation*}
\Delta_j = (s_{j-1} \a, s_j \a],
\ j = 1, 2, \dots, L.
\end{equation*}
Define step functions $q_{n,1}, q_{n, 2}$:
\begin{gather*}
q_{k, 1}(t) = q_k(s_{j-1} \a)
= \frac{1}{g(s_{j-1}+ \b)^{2\g}},\ t\in \Delta_j, \\
q_{k, 2}(t) = q_k(s_j \a)
= \frac{1}{g(s_j+ \b)^{2\g}},\ t\in \Delta_j.
\end{gather*}
Here we have denoted $\b = \b_k$.
Further proof is for the case $\sign = -1$ only.
The other case
is done in the same way.
Clearly, $q_{k, 2}\le q_k \le q_{k, 1}$, and hence
the counting function
of the operator $A_0 - gq_k$
with the Dirichlet conditions at the ends of the interval
$\Delta$ satisfies the two-sided estimate
\begin{equation*}
N(\l; -g q_{k, 2}, \Delta)
\le N(\l; -gq_k, \Delta)\le
N(\l; -gq_{k, 1}, \Delta).
\end{equation*}
Let us find the asymptotics of the r.h.s.
We are going to use the decoupling principle again:
\begin{equation}\label{above:eq}
N(\l; -gq_{k, 1}, \Delta)\le \sum_{j=1}^L
N(\l; -gq_{k, 1}, \Delta_j) + 2(L-1).
\end{equation}
Now, using Theorem \ref{density:thm}, we get for each $j$
\begin{align*}
\a^{-1} N(\l; -g q_{k, 1}, \Delta_j )
= &\ (s_j - s_{j-1}) |\Delta_j|^{-1}
N\bigl(\l + (s_{j-1} + \b)^{-2\g} ; 0, \Delta_j\bigr)\\
\le  &\ (s_j - s_{j-1})\rho\bigl(
\l + (s_{j-1} + \b)^{-2\g}
\bigr)\\
&\  + C\bigl(1+ \sqrt{|\l|
+(s_{j-1} + \b)^{-2\g}}\bigr)\a^{-1},
\end{align*}
with a universal constant $C$. Since $\b \ge \d$, we obtain
from \eqref{above:eq} that
\begin{multline*}
\a^{-1} N(\l; -g q_{k,1}, \Delta)
- \sum_j (s_j - s_{j-1})
\rho\bigl(
\l + (s_{j-1} + \b)^{-2\g}
\bigr)\\
\le(L + C(\l)+ C\d^{-\g}) \a^{-1}.
\end{multline*}
This can be rewritten as
\begin{multline}\label{for:eq}
\a^{-1} N(\l; -g q_{k,1}, \Delta)
- \sum_j \int_{\Delta_j}
\rho\bigl(
\l + (s_{j-1} + \b)^{-2\g}
\bigr) ds\\
\le(L + C(\l)+ C\d^{-\g}) \a^{-1}.
\end{multline}
To replace the sum in the l.h.s.
by the integral,
we use the H\"older property \eqref{holder1:eq}
with
\begin{equation*}
\l_1 = \l + (s_{j-1} + \b)^{-2\g}, \
\l_2 = \l + (t + \b)^{-2\g}.
\end{equation*}
Observe that
\begin{align*}
|\l_1 - \l_2| = &\ |(s_{j-1} + \b)^{-2\g} - (t+\b)^{-2\g}|\\
\le &\ 2\g\d^{-2\g-1}|s_{j-1} - t|
\le 2\g R\d^{-2\g-1}L^{-1},\ \forall t\in\Delta_j.
\end{align*}
Now we infer from \eqref{holder1:eq} that
\begin{equation*}
\biggl|\rho\bigl(
\l + (s_{j-1} + \b)^{-2\g}
\bigr) - \rho\bigl(
\l + (t + \b)^{-2\g}
\bigr)\biggr|
\le  C \d^{-\g - 1/2} R^{1/2} L^{-1/2},
\ t\in \Delta_j.
\end{equation*}
Substituting this estimate into \eqref{for:eq},
we get
\begin{gather*}
\a^{-1} N(\l; -g q_{k, 1}, \Delta)
- \int_0^R
\rho\bigl(
\l + (s + \b)^{-2\g}
\bigr) ds \le C E(\a, \d,  L; R), \\
E(\a, \d, L; R) =  (L + 1 + \d^{-\g}) \a^{-1}
+ \d^{- \g - 1/2} R^{3/2} L^{-1/2}.
\end{gather*}
Arguing similarly, we arrive at the analogous
lower bound for $N(\l; -g q_{k, 2}, \Delta)$.
Consequently,
\begin{equation*}
\biggl|\a^{-1} N(\l; -g q_k, \Delta)
- \int_0^R
\rho\bigl(
\l + (s + \b)^{-2\g}
\bigr) ds\biggr|\le CE(\a, \d, L; R)
\end{equation*}
In view of Theorem \ref{density:thm} we also have
\begin{equation*}
\bigl|\a^{-1} N(\l; 0, \Delta)
- \int_0^R \rho(\l) ds\bigr|
\le C \a^{-1}.
\end{equation*}
By \eqref{w:9},
in combination with the previous estimate this gives
\begin{equation*}
\biggl|\a^{-1} M(\l; -g q_k, \Delta)
- F_-^R(\b_k, \l)
\biggr|
\le CE(\a, \d, L; R)
\end{equation*}
The parameter $L$ can be chosen so as
to insure that $E\to 0$ as $\a \to \infty$.
Indeed, in view of \eqref{small:eq} $\a\d^{2\g+1}\to \infty$
as $\a\to\infty$.
Therefore, defining
\begin{equation*}
L = [\a^{1/2} \d^{-\g - 1/2}],
\end{equation*}
we guarantee that
\begin{equation*}
\d^{-2\g-1} L^{-1} \sim \a^{-1/2}\d^{-\g-1/2}\to 0,
\ \ L\a^{-1}\sim \a^{-1/2} \d^{-\g-1/2}\to 0,
\ \a\to\infty.
\end{equation*}
Taking $R$ to infinity and referring to
Lemma \ref{f:lem}(iv),
we obtain \eqref{final:eq}, thus completing the proof
for $q(t) = Q(t)$.

It remains to include the potentials satisfying \eqref{asympt:eq}.
To this end note that under the condition \eqref{asympt:eq},
for any $\vare >0$
\begin{equation*}
Q_k(t)(1-\vare) \le q_k(t)\le Q_k(t)(1+\vare),
\end{equation*}
if $k$ is sufficiently large.
Thus, using the monotonicity of $M(\l; V)$ in $V$
(see Sect. \ref{tree4:sect}) and the asymptotics
\eqref{final:eq} for $q = Q$ we easily deduce
\eqref{final:eq} for the general case.
\end{proof}

In conclusion note that we shall not need
the asymptotics \eqref{final1:eq} in what follows.

\section{Asymptotics of $\widetilde M(\l, \BA_{\pm gq})$:
proof of Theorems
\ref{main:thm}, \ref{g2:thm}, \ref{main1:thm},
\ref{exp+:thm}, \ref{exp-:thm}}
\label{tree8:sect}

Throughout this section we assume that $\l, \l_1, \l_2$
are $2$-admissible and satisfy \eqref{away:eq}.

\subsection{Proof of Theorems \ref{main:thm}
and \ref{main1:thm}}
Recall that in Theorem \ref{main:thm}
we assume that either $\sign = -1,\ \g <2$, or $\sign = +1$ and
$\g>0$ is arbitrary. In Theorem \ref{main1:thm}
$\g >0$ is arbitrary, but if $\sign = -1$ and $\g \ge 2$, then
the potential $q$ satisfies Condition \ref{q:cond}.

\underline{Step I.}\  To begin with we
show that ``small'' or ``large''
values of $k$ do not contribute to the
sums $\widetilde M$ and $\widetilde N$.

Suppose that $k\ge A\a$ with some fixed $A > 0$.
Then for $\g \le 2$  and large $A$
the perturbation $gq_k\le CA^{-2\g}$
is small and therefore $M(\l, \pm gq_k) = 0$, since $\l$ is
strictly inside the gap. For $\g >2$,
by \eqref{known:eq} we have
\begin{equation*}
\sum_{k\ge A\a} M(\l;  \pm gq_k)\le Cg^{1/2}(A\a)^{2-\g}
= C'A^{2-\g}\a^2.
\end{equation*}
By \eqref{t:eq} a similar bound holds for the
sum of the functions
$N(\l_1, \l_2; \pm gq_k)$ for all $\g > 0$.
These calculations again show that
$k\ge A\a$ do not contribute as $A$ grows.

Suppose that $k \le \d \a$.
If $\g < 2$ and $\sign = -1$, then
the bounds \eqref{known:eq} (for $\g>1$) and
\eqref{estpm:eq} (for $\g\le1$) ensure that for $\d>0$
\begin{equation}\label{novoe:eq}
\sum_{k\le \d\a} M(\l; -gq_k)\le
\begin{cases}
C\d^{2-\g} \a^2,\ 1<\g < 2;\\[0.2cm]
C\d\a^2\ln\a,\ \g = 1;\\[0.2cm]
C\d\a^2,\ \g < 1.
\end{cases}
\end{equation}
This means that the share of this
sum becomes small when $\d\to 0$.
For $\g\not=1$ we can take $\d$ to be arbitrarily small
constant, \textsl{independent of} $\a$. With
$\g=1$ we must be more careful. Since we want to obtain
the asymptotics of order $\a^2$, we should ``kill'' the ``$\ln$''
term in the estimate by choosing $\d$ to be dependent on $\a$,
but in a very mild way: $\d = \a^{-\eta}$ with a parameter
$\eta<(1+2\g)^{-1}$, so that the condition \eqref{small:eq}
from Lemma \ref{individ:lem} is satisfied.
For $\sign = +1$
the estimate \eqref{est+:eq} yields:
\begin{equation*}
\sum_{k\le \d\a}M(\l; gq_k)\le C_\g\d \a^2,\ \forall \g >0.
\end{equation*}
As in the case $\sign=-1$ and $\g\neq1$, it is possible to take
$\d$ to be an arbitrarily small constant.
However, for the sake of uniformity,
we take $\d=\a^{-\eta}$, for both signs
$\sign=\pm1$.
Consequently,
\begin{equation*}
\sum_{k\le \d\a}M(\l; \pm gq_k)=o(\a^2)
\end{equation*}
and the condition \eqref{small:eq} is satisfied.
In the case of the function $\widetilde N$
the estimate \eqref{est:eq} guarantees that
\begin{equation}\label{novoe1:eq}
\sum_{k\le \d\a}N(\l_1, \l_2; \sign gq_k)
\le C_\g\d \a^2 = o(\a^2),
\end{equation}
for $\sign = +1$ and all $\g >0$. If
$\sign = -1,\  \g < 2$, then the same
bound follows from \eqref{novoe:eq} and \eqref{t:eq}.
In the case $\sign = -1,\  \g \ge 2$ the estimate \eqref{novoe1:eq}
is a direct consequence of \eqref{n:eq}.
Thus, it remains to study the
sums \eqref{w:4}, \eqref{w:5} only over the numbers
\begin{equation*}
[\d \a]\le k\le [A \a],
\end{equation*}
with $\d = \a^{-\eta},\ \eta < (1+2\g)^{-1}$, and a
fixed $A \ge \sup_{\a}\d$.

\underline{Step 2.}
We use the notation
\eqref{ab:eq}.
Estimate using \eqref{sumf:eq}:
\begin{align*}
\limsup \biggl|\a^{-2}\sum_{k=[\d\a]}^{[A\a]}&\
M(\l; \pm g q_k)  -  \int_{\d}^A F_\pm(\s,\l) d\s
\biggr|\\
\le &\ \limsup \a^{-1}\sum_{k=[\d\a]}^{[A\a]}
\biggl| \a^{-1} M(\l; \pm gq_k)
- F_\pm(\b_k, \l)\biggr| \\
\le &\  A\limsup \max_{[\d\a]\le k\le [A\a]}
\biggl| \a^{-1} M(\l; \pm gq_k)
- F_\pm(\b_k, \l)\biggr|.
\end{align*}
The r.h.s. tends to zero by Lemma \ref{individ:lem}.
Consequently
\begin{equation*}
\lim \biggl[\a^{-2}\sum_{k=[\d\a]}^{[A\a]}
M(\l; \pm g q_k)  -
 \int_{\d}^A F_\pm(\s,\l) d\s\biggr] = 0,\ g\to\infty.
\end{equation*}
By \eqref{t:eq} and \eqref{gpm:eq}
this equality implies that
\begin{equation*}
\lim \biggl[\a^{-2}\sum_{k=[\d\a]}^{[A\a]}
N(\l_1, \l_2; \pm g q_k) -
 \int_{\d}^A G_\pm(\s,\l_1, \l_2) d\s\biggr] = 0,\ g\to\infty.
\end{equation*}
Referring to Step I of the proof and Lemma \ref{sumf:lem}(ii),
we can now replace the lower and upper limits of
summation and integration by $0$ and $\infty$ respectively.
This completes
the proof of Theorems \ref{main:thm}, \ref{main1:thm}.
\qed

\subsection{Proof of Theorem \ref{g2:thm}}
By \eqref{asympt:eq},
\begin{equation*}
\int_0^\infty\sqrt{q_k(t)}dt
= \int_0^\infty\frac{1+o(1)}{(1+k+t)^2}dt=
\frac{1+\e_k}{k+1},
\end{equation*}
where $\e_k\to 0$ as $k\to\infty$.
From here we find by virtue of Theorem \ref{weyl1:thm}
that the function $M((\l)$ has the following asymptotics:
\begin{equation}\label{g2:eq}
\biggl|M(\l;-gq_k)-\frac{g^{1/2}}{\pi(k+1)}\biggr|
\le Cg^{1/4} + C'g^{1/2}\frac{\e_k}{k+1},
\end{equation}
uniformly in $k\ge 0$. Observe that
the components with numbers $k\ge A g^{1/4}$ do
not contribute if $A$ is sufficiently
large, since $|gq_k(t)|\le
CA^{-4}$. Let us turn to the remaining terms:
\begin{multline*}
\biggl|\widetilde M(\l,  \BA_{-gq})
- \frac{g^{1/2}\ln g}{4\pi}
\biggr|
\le \sum_{k\le Ag^{1/4}}\bigl|M(\l; -gq_k)
- g^{1/2} (k+1)^{-1} \pi^{-1}\bigr| \\
+ g^{1/2}\pi^{-1}\biggl|\sum_{k\le Ag^{1/4}}
(1+k)^{-1} - \ln g/4\biggr|.
\end{multline*}
The second term in the r.h.s. is of order $O(g^{1/2})$ in
view of the known formula for
the partial sum of the harmonic
series.
By \eqref{g2:eq} the first term
is bounded by
\begin{equation*}
CAg^{1/2} + C'g^{1/2} \sum_{k\le Ag^{1/4}} \frac{\e_k}{1+k}.
\end{equation*}
Since $\e_k\to 0$ as $k \to \infty$, this quantity is
of order $o(g^{1/2}\ln g)$.
\qed

\subsection{Proof of Theorems
\ref{exp+:thm}, \ref{exp-:thm}}
Rewrite the sum $\widetilde M$ in the form
\begin{equation*}
\widetilde M(\l, \BA_{gq}) = \sum_k M\bigl(\l; ge^{-2\vark k}
\hat q^{(k)}\bigr),\ \hat q^{(k)} =  e^{2\vark k}q_k.
\end{equation*}
Since $\hat q^{(k)}, k \ge 0,$
satisfies the bound \eqref{twoside:eq}
for all $t\ge R_0$, from Theorem \ref{plus:thm}
we obtain that
\begin{equation*}
\bigl|M\bigl(\l; ge^{-2\vark k} \hat q^{(k)}\bigr)
- \rho(\l) (\a - k) \bigr|\le C,\ \forall k\le \a-1.
\end{equation*}
Consequently,
\begin{align*}
\sum_{k\le \a-1} M(\l; ge^{-2\vark k} \hat q^{(k)})
= &\ \rho(\l) \sum_{k\le \a-1} (\a-k) + O(\a)\\
= &\ \frac{1}{2}\rho(\l)\a^2 + O(\a).
\end{align*}
On the other hand, by \eqref{known1:eq}
\begin{equation*}
M\bigl(\l; \pm g q_k\bigr)\le Cg^{1/2}
e^{-\vark k} = C  e^{\vark(\a-k)},
\end{equation*}
so that
\begin{equation}\label{largen:eq}
\sum_{k > \a-1} M(\l; \pm g q_k) \le C.
\end{equation}
The asymptotics \eqref{tildemexp:eq}  follows.

The estimates \eqref{tildenexp+:eq} and \eqref{tildenexp-:eq}
are proved in
the same way. By \eqref{nexp:eq} and
\eqref{largen:eq}, \eqref{t:eq}
\begin{equation*}
\widetilde N(\l_1, \l_2; \BA_{gq})
\le \sum_{k\le \a - 1} C + C'\le C\a.
\end{equation*}
Furthermore, as $\hat q^{(k)}$ satisfies Condition
\ref{q:cond},
by \eqref{n:eq} and \eqref{largen:eq}, \eqref{t:eq}
\begin{align*}
\widetilde N(\l_1, \l_2; \BA_{-gq})
\le &\ \sum_{k\le \a-1} N(\l_1, \l_2; -ge^{-2\vark k} \hat q_k)
+ C'\\
\le &\ \sum_{k\le \a - 1}C(\a-k) + C'\le C''\a^2.
\end{align*}
\qed

\section{Asymptotics of $M(\l, \BA_{\pm gq})$
and $N(\l, \BA_{\pm gq})$}
\label{tree9:sect}

Here we turn to the study of the sums
\eqref{w:2} and \eqref{w:3}. As before, to ensure that
they are finite we assume that
$\l, \l_1, \l_2$ satisfy \eqref{away:eq}
with the same closed interval $I$.
Due to the presence of exponential terms in the sums,
their study is more complicated than that of
$\widetilde M$, $\widetilde N$, and hence the asymptotic
formulae are less explicit.
Another feature is that
for the exponential and power-like potentials the results
are qualitatively different.

\subsection{Exponential potentials} In this subsection
we always (except for Theorem \ref{critic:thm}) assume that
\begin{equation}\label{exact:eq}
q(t) = Q(t) = e^{-2\vark t}.
\end{equation}
This assumption allows one to obtain
asymptotic formulae based on the
``self-similarity'' property
of the function $e^{-2\vark t}$.
Introduce the notation
\begin{equation*}
\ln b = \b >0.
\end{equation*}

\begin{thm}\label{multiple:thm}
Assume \eqref{exact:eq}.
Then the following two statements hold:
\begin{itemize}
\item[(i)]
Let $\vark > 0$ be arbitrary.
Then there exist functions $\varphi_\pm$ that are $2\vark$-periodic,
bounded and separated from zero, such that
\begin{equation}\label{phi:eq}
\lim \bigl[g^{-\frac{\b}{2\vark}}N(\l_1, \l_2; \BA_{\pm gq}) -
\varphi_\pm (\ln g)\bigr]  = 0,\ g\to\infty.
\end{equation}
\item[(ii)]
Suppose that $\vark > 0$ is arbitrary if $\sign  = +1$
and $\vark < \b$ if $\sign = -1$. Then
there exist two functions
$\psi_\pm$ that are $2\vark$-periodic,
bounded and separated from zero, such that
\begin{equation}\label{psi:eq}
\lim \bigl[g^{-\frac{\b}{2\vark}}M(\l, \BA_{\pm gq} ) -
\psi_\pm(\ln g)\bigr] = 0,\ g\to\infty.
\end{equation}
\end{itemize}
\end{thm}

We precede the proof with an elementary but convenient lemma:

\begin{lem}\label{renew:lem}
Let $n(t), t\in \R$ be a bounded function
such that
\begin{equation}\label{ncond:eq}
\begin{cases}
n(t) = 0\ \  \textup{ for all} \ \
t\le t_0 \ \ \textup{with some}\ \  t_0 > 0,\\
n(t)\le C t^{\frac{\b}{2\vark} - \e},
t\ge t_0, \ \textup{for some}\ \ \e >0.
\end{cases}
\end{equation}
Then for the function
\begin{equation}\label{ndef:eq}
N(g) = n(g) + (1-b^{-1})\sum_{k\ge 1} e^{\b k} n(g e^{-2\vark k})
\end{equation}
there exists a function $\phi$ which is $2\vark$-periodic,
bounded and separated from zero, such that
\begin{equation*}
\lim \bigl[g^{-\frac{\b}{2\vark}}N(g) -
\phi(\ln g)\bigr] = 0,\ g\to\infty.
\end{equation*}
\end{lem}

\begin{proof}
The sum in the r.h.s. of \eqref{ndef:eq}
is finite, since for sufficiently
large $k$ we have $ge^{-2\vark k} \le t_0$.
Denote
\begin{equation*}
\xi(g) = g^{-\frac{\b}{2\vark}}n(g),\ \
\Xi(g) = g^{-\frac{\b}{2\vark}}N(g).
\end{equation*}
Then \eqref{ndef:eq} yields
\begin{equation*}
\Xi(g) = \xi(g) + (1-b^{-1}) \sum_{k\ge 1} \xi(g e^{-2\vark k}),
\end{equation*}
which in its turn implies that
\begin{equation*}
\Xi(g) - \Xi (ge^{-2\vark}) =
\xi(g) - b^{-1} \xi(g e^{-2\vark}).
\end{equation*}
Using the notation \eqref{alpha:eq}
and introducing new functions
$F(\a) = \Xi(g)$, $f(\a) = \xi(g)$, we arrive at the equation
\begin{equation*}
F(\a) - F(\a-1) = f(\a) - b^{-1} f(\a-1).
\end{equation*}
Since $n(t)$ satisfies \eqref{ncond:eq},
$f(\a) = 0$ for $\a\le \a_0 = (2\vark)^{-1}\ln t_0$, and
$f(\a)\le Ce^{-2\vark \e \a}$, $\a\ge \a_0$.
Therefore all the conditions of the Renewal Theorem
are satisfied
(see \cite{F}, Chapter XI.1,
or a modern exposition in
\cite{LV}), which guarantees the
existence of a $1$-periodic
function $\widetilde\phi$ which is bounded
and separated from zero, such that
\begin{equation*}
F(\a) = \widetilde\phi(\a) + o(1),\  \a\to\infty.
\end{equation*}
which leads to \eqref{phi:eq} after substitution
$\phi(\a) = \widetilde\phi\bigl(\a/(2\vark)\bigr)$.
\end{proof}

\begin{proof}[Proof of Theorem \ref{multiple:thm}]
(i) The proof is done simultaneously
for both signs $\sign = \pm 1$. We use Lemma
\ref{renew:lem} with $n(g) = N(\l_1, \l_2; \pm gq)$.
Since $N(\l_1, \l_2; \pm gq_k) = n(ge^{-2\vark k})$,
by \eqref{w:3} the function $N(g)$ in the r.h.s.
of \eqref{ndef:eq} coincides with $N(\l_1, \l_2; \BA_{\pm gq})$.
By \eqref{known1:eq} and \eqref{t:eq} $n(t) = 0, \ t\le t_0$
for a sufficiently small $t_0 >0$.
Moreover, by
\eqref{n:eq} or \eqref{est:eq}, the second condition in
\eqref{ncond:eq} is also fulfilled for any
$\e < \b(2\vark)^{-1}$. Thus the required
asymptotics \eqref{phi:eq} follows from Lemma \ref{renew:lem}.

(ii) The cases $\sign = +1$ and
$\sign = -1$ are treated separately.

Let first $\sign = +1$.
Denote now $n(g) = M(\l; gq)$.
Then, similarly to the first part of the proof,
the total counting function \eqref{w:2}
coincides with \eqref{ndef:eq}.
By \eqref{known1:eq} and
\eqref{est+:eq} the function $n$
satisfies \eqref{ncond:eq} for any $\e < \b(2\vark)^{-1}$.
Thus Lemma \ref{renew:lem}
guarantees the asymptotics \eqref{psi:eq} for $\sign = +1$.

In the case $\sign = -1$, $\vark < \b$, the first
condition in \eqref{ncond:eq} is satisfied for
$n(g) = M(\l; -gq)$ in view of \eqref{known1:eq}.
Besides, \eqref{known1:eq} ensures also that
the second condition is satisfied with
$\e = \b(2\vark)^{-1} -1/2 >0$.
Again, Lemma  \ref{renew:lem}
leads to \eqref{psi:eq} for $\sign = -1$.
\end{proof}

For $\sign = -1$ the cases $\vark < \b$ and
$\vark > \b$ are described by Theorem
\ref{multiple:thm} and Lemma \ref{expo:lem}
respectively. Let us handle the critical case
$\b = \vark$. We emphasise that this is the only
asymptotic formula in this subsection which does not
require the exact equality $q(t) = e^{-2\vark t}$.

\begin{thm}\label{critic:thm}
Suppose that
\begin{equation*}
q(t) = Q(t)(1+o(1)),\ t\to\infty,\ \ Q(t) = e^{-2\vark t},
\end{equation*}
with $\vark = \b$, and
that $q$ satisfies Condition
\ref{q:cond}.
Then
\begin{equation*}
\lim\dfrac{M(\l; \BA_{-gq})}{g^{1/2} \ln g} =
\frac{1-b^{-1}}{2\pi\vark^2}, \ g\to\infty.
\end{equation*}
\end{thm}

\begin{proof}
By \eqref{known1:eq}
$M(\l; -gq_k) = 0$ for all $k\ge \a + A$ with
a sufficiently large $A$, and
\begin{equation*}
\sum_{k\le A} b^k M(\l; -gq_k)
+ \sum_{\a-A\le k\le \a + A} b^k M(\l; -gq_k)
\le CA g^{1/2}.
\end{equation*}
Consequently
\begin{equation}\label{critic:eq}
\lim\dfrac{M(\l; \BA_{-gq})}{g^{1/2} \ln g} =
\lim
\frac{1-b^{-1}}{g^{1/2} \ln g}\sum_{A< k< \a - A}
e^{\vark k} M(\l; - gq_k), \ g\to\infty,
\end{equation}
if the limit in the r.h.s. exists.
For $k\in (A, \a-A)$ apply Theorem \ref{weyl1:thm} and the
relation \eqref{shift:eq}
to obtain the asymptotics
\begin{align*}
M(\l; -gq_k) = &\ \frac{g^{1/2}}{\pi}
\int_0^\infty \sqrt{q(k+t)} dt
+ O(\a)\\
 = &\ g^{1/2} e^{-\vark k}\bigl((\pi\vark)^{-1} + o_A(1)\bigr)
 + O(\a),
\end{align*}
where $o_A(1)\to 0$ as $A\to\infty$ uniformly in $k, g$.
Therefore the sum
in the r.h.s. of \eqref{critic:eq} equals
\begin{gather*}
g^{1/2}\sum_{A < k < \a - A}
\frac{1}{\pi\vark} + g^{1/2}\a o_A(1)
+ O(\a)\sum_{A < k <\a-A} e^{\vark k}\\
=  \frac{g^{1/2}\a}{\pi\vark} + g^{1/2} O(A)
+  g^{1/2}\a o_A(1)
+ e^{\vark\a} e^{-\vark A} O(\a).
\end{gather*}
Since $\a = \ln g/(2\vark)$,
now it follows from \eqref{critic:eq} that
\begin{equation*}
\limsup_{g\to\infty}\biggl|
\dfrac{M(\l; \BA_{-gq})}{g^{1/2} \ln g}
- \frac{1}{2\pi\vark^2}\biggr| = o_A(1) + O(e^{-\vark A}).
\end{equation*}
Since $A$ is arbitrary, the required result follows.
\end{proof}

\subsection{Power-like potentials}
For power-like potentials
the asymptotic formulae that we obtain,
are less informative since
they are established for $\ln M$ and $\ln N$.
For the sake of illustration we consider
here only $M(\l; \BA_{\pm g q})$.
The corresponding asymptotics of
$N(\l; \BA_{\pm gq})$ can be easily derived using the
same argument as well.
For simplicity we assume that
$q(t) = Q(t) = (1+t)^{-2\g}$.
For more general power-like potentials
the results follow by
monotonicity of $M$
with respect to the potential.
Recall that $I$ denotes the interval defined in
\eqref{away:eq}. We also use the notation
$\dc_\pm = |\l-\l_{\mp}|$,
$\s_\pm = \dc_{\pm}^{-\frac{1}{2\g}}$ introduced in
Lemma \ref{f:lem}(iii).

\begin{thm} Let $I$ be a closed interval
defined in \eqref{away:eq}, which is strictly inside
the gap $(\l_-, \l_+)$.
Let $q(t) = (1+t)^{-2\g},\  \g >0$.
Then for any $\l\in I$
\begin{equation*}
\lim \a^{-1} \ln M(\l; \BA_{\pm gq}) = \b \dc_{\pm}^{-\frac{1}{2\g}},\
g\to\infty,
\end{equation*}
where $\b = \ln b$.
\end{thm}

\begin{proof}
\textsl{Upper bound.}
A straightforward perturbation
argument ensures that $M(\l; \pm gq_k) = 0$
if $g q_k(t)\le \dc_{\pm},\ \forall t >0$, i.e. for
\begin{equation*}
k\ge K_1 = K_1(g) = (\dc^{-1} g )^{\frac{1}{2\g}}
= \a \dc_{\pm}^{-\frac{1}{2\g}}.
\end{equation*}
It follows from
Lemmas \ref{power:lem} and \ref{est:lem} that
$M(\l; \pm gq_k)\le Cg^\om$ with some $\om = \om(\g)>0$.
Substituting this bound in \eqref{w:2}, gives
\begin{align*}
M(\l; \BA_{\pm gq}) = &\ M(\l; \pm qq)
+ (1-b^{-1})\sum_{k\in\mathbb N}b^k M(\l; \pm gq_k)\\
\le &\ C g^{\om}\sum_{k = 0}^{K_1} b^k\le C' g^\om b^{K_1},
\end{align*}
with a constant $C'$ depending only on $b$, $\g$.

\textsl{Lower bound.}
For the lower bound we drop all but one term
from the sum \eqref{w:2}: for any $\e\in(0, 1)$
we have
\begin{equation*}
M(\l; \BA_{\pm gq})
\ge  (1-b^{-1})b^k M(\l; \pm gq_k),\ k = [(1-\e) K_1].
\end{equation*}
By Lemmas \ref{individ:lem} and \ref{f:lem}(iii),
\begin{equation*}
M(\l; \pm gq_k) \ge c\a F_{\pm}(\l; \b_k)
\ge c'\a\bigl(\dc_{\pm}^{-\frac{1}{2\g}} - \b_k\bigr)^{\frac{3}{2}},
\end{equation*}
for sufficiently large $\a$. Since $\b_k = (k+1)\a^{-1}$,
we see that
the r.h.s. is bounded from below by $c'' \a\e^{3/2}$
with a constant depending only on $\dc_{\pm}$.
Consequently,
\begin{equation*}
M(\l; \BA_{\pm gq})\ge c \e^{\frac{3}{2}}\a b^{K_1(1-\e)},
\end{equation*}
with a constant depending on $b, \l$ and $\g$.
Since $\e > 0$ is arbitrary,
in combination with the upper bound,
this gives the required asymptotics.
\end{proof}

\section{Acknowledgements}

The second author (M.S.) was partly supported by the
Minerva center for non-linear physics and by the
Israel Science Fundation,
and partly by the EPSRC grant GR/N 37193/01.

This paper was essentially completed when M.S. was visiting
King's College, London in April-May 2001.
The authors are grateful to Yu. Safarov for
valuable discussions.

\bibliographystyle{amsplain}

\providecommand{\bysame} {\leavevmode\hbox
to3em{\hrulefill}\thinspace}

\end{document}